\documentclass[14pt, draft]{article}
\usepackage[cp1251]{inputenc}
\usepackage[ukrainian, russian, english]{babel}

\usepackage{ amsfonts, amssymb}
\usepackage{amsmath}
\usepackage{amsthm}

\begin{document}
\selectlanguage{ukrainian} \thispagestyle{empty}
 \pagestyle{myheadings}              

УДК 517.51 \vskip 3mm

\bf Сердюк А.С., Мусієнко А.П.  \rm (Інститут математики НАН України, Київ)

\bf A.S. Serdyuk, A.P. Musienko\footnote{Робота виконана за часткової підтримки Державного фонду фундаментальних досліджень України (проект \mbox{№ GP/Ф36/068)}} \rm

\vskip 5mm

\centerline{\bf Нерівності типу Лебега для сум Валле Пуссена на множинах} \centerline {\bf  аналітичних та цілих функцій}
\vskip 4mm
\centerline  {\bf Lebesgue-type inequalities  for de la  Vall\'{e}e Poussin sums on the sets} \centerline {\bf of analytic and entire functions}
\vskip4mm

\vskip 5mm

\rm For the functions from sets $C_\beta^\psi C$ and $C_\beta^\psi L_s, \ 1\leq s\leq\infty$, generated by sequences $\psi(k)>0$ satisfying the condition d'Alembert $\mathop {\rm \lim}\limits_{k\rightarrow\infty}\frac{\psi(k+1)}{\psi(k)}=q, \ q\in[0,1)$, asymptotically unimprovable estimates for deviations  of de la Vall\'{e}e Poussin sums in the uniform metric, which are represented in terms of values of the best approximations of $(\psi,\beta)$-differentiable functions of this sort by trigonometric polynomials in the metrics $L_s$ are obtained. Proved that received estimates are unimprovable on some important functional subsets.
\vskip 5mm
\rm Для функцій з множин $C_\beta^\psi C$ та $C_\beta^\psi L_s, \ 1\leq s\leq\infty$, породжуваних послідовностями $\psi(k)>0$, що задовольняють умову Даламбера ${\mathop {\rm \lim}\limits_{k\rightarrow\infty}\frac{\psi(k+1)}{\psi(k)}=q, \ q\in[0,1)}$, одержано асимптотично непокращувані оцінки відхилень в рівномірній метриці сум Валле Пуссена, які виражаються через значення найкращих наближень $(\psi,\beta)$-похідних таких функцій тригонометричними поліномами в метриках просторів $L_s$. Доведено, що одержані оцінки залишаються непокращуваними  на деяких важливих функціональних підмножинах з $C_\beta^\psi C$ та $C_\beta^\psi L_s$.
\vskip 5mm

\rm  Для функций из множеств $C_\beta^\psi C$ и $C_\beta^\psi L_s, \ 1\leq s\leq\infty$, порождаемых последовательностями $\psi(k)>0$, удовлетворяющих условие Даламбера $\mathop {\rm \lim}\limits_{k\rightarrow\infty}\frac{\psi(k+1)}{\psi(k)}=q, \ q\in[0,1)$, получены асимптотически неулучшаемые оценки отклонений в равномерной метрике сумм Валле Пуссена, выражающиеся через значения наилучших приближений $(\psi,\beta)$-производных этих функций тригонометрическими полиномами в метриках пространств $L_s$. Доказана,  неулучшаемость  полученных оценок на некоторых важных функциональных подмножествах из $C_\beta^\psi C$ и $C_\beta^\psi L_s$.

\newpage

Нехай $C$ --- простір  $2\pi$-періодичних неперервних функцій $f$ з
нормою ${\|f\|_C=\mathop {\rm max}\limits_{t}|f(t)|;}$
$L_s$ --- простір  $2\pi$-періодичних сумовних  функцій з нормою
$$
\|f\|_s=\|f\|_{L_s}= \left\{\begin{array}{ll}
\Big(\int\limits_0^{2\pi}|f(t)|^sdt\Big)^\frac{1}{s},& \ \ 1\leq s<\infty,\\
\mathop{\rm ess\,sup}\limits_{0\leq t\leq 2\pi}|f(t)|,&  \ \  s=\infty.
\end{array}\right.
$$

В даній роботі розглядаються множини $L^\psi_\beta\mathfrak{N}$ і $C^\psi_\beta\mathfrak{N}$ періодичних функцій, які було введено О.І. Степанцем [\ref{Step monog 1987}, \ref{Step monog 2002(1)}] наступним чином. Нехай $f$ --- $2\pi$-періодична, сумовна  функція (${f\in L_1}$) і
$$
\frac{a_0}{2}+\sum\limits_{k=1}^\infty(a_k\cos kx+b_k\sin kx)\equiv\sum\limits_{k=0}^\infty A_k(f;x)
$$
--- її ряд Фур'є. Нехай, далі, $\psi(k)$ --- довільна послідовність дійсних чисел і $\beta$~--- фіксоване дійсне число.

Якщо ряд
$$
\sum\limits_{k=1}^\infty\frac{1}{\psi(k)}\bigg(a_k\cos\big(kx+\frac{\beta\pi}{2}\big)+b_k\sin\big(kx+\frac{\beta\pi}{2}\big)\bigg)
$$
є рядом Фур'є деякої сумовної функції $\varphi$, то цю функцію називають $(\psi,\beta)$-похідною функції $f$ і позначають через $f_\beta^\psi$ [\ref{Step monog 1987}, c. 25]. Множину всіх функцій $f$, які задовольняють таку умову, позначають через $L_\beta^\psi$. Якщо $f\in L_\beta^\psi$ і в той же час  $f_\beta^\psi\in\mathfrak{N}$, де $\mathfrak{N}$ --- деяка підмножина з ${L_1^0=\{\varphi\in L_1:\int\limits_{-\pi}^{\pi }\varphi(t) dt=0\}}$, то записують $f\in L_\beta^\psi\mathfrak{N}$. Якщо $F_\beta^\psi=f$, то функцію $F$ називають $(\psi,\beta)$-інтегралом функції $f$, при цьому записують $F(x)={\cal J}_\beta^\psi(f;x)$. Покладемо $L_\beta^\psi\bigcap C=C_\beta^\psi$, $L_\beta^\psi\mathfrak{N}\bigcap C=C_\beta^\psi\mathfrak{N}$. Надалі в ролі $\mathfrak{N}$ виступатимуть простори $C$ чи ${L_s, \  1\leq s\leq\infty}$, або ж деякі їх підмножини.

В рамках роботи будемо вважати, що послідовність $\psi(k)$, яка породжує класи $L_\beta^\psi\mathfrak{N}$ і $C_\beta^\psi\mathfrak{N}$ додатна і  задовольняє умову ${\cal D}_q,\,\,\, q\in[0,1)$, яка полягає у виконанні наступної рівності:
\begin{equation}\label{S4.0.4}
\mathop {\rm \lim}\limits_{k\rightarrow\infty}\frac{\psi(k+1)}{\psi(k)}=q.
\end{equation}
В такому разі будемо використовувати запис: $\psi\in{\cal D}_q$.  Множини $C^\psi_\beta\mathfrak{N}$, $\psi\in{\cal D}_q$, $q\in(0,1)$  складаються з $2\pi$-періодичних функцій $f(x)$, які допускають регулярне продовження у смугу $|\mbox{\rm Im} z|\leq \ln\frac{1}{q}$ комплексної площини. У випадку, коли $\psi\in {\cal D}_0,$ множини $C^{\psi}_\beta\mathfrak{N}$ складаються з функцій, регулярних в усій комплексній площині, тобто з цілих функцій.

Відомо [\ref{Step monog 2002(1)}, c. 136], що класи $L_\beta^\psi\mathfrak{N}, \ \psi\in {\cal D}_q, \ q\in[0,1)$ складаються з функцій, які майже скрізь можна зобразити у вигляді згортки
\begin{equation}\label{0.3}
f(x)=\frac{a_0}{2}+\frac{1}{\pi}\int\limits_{-\pi}^{\pi}\varphi(x-t)\Psi_\beta(t)dt, \ \varphi\in\mathfrak{N}, \ \varphi\perp1
\end{equation}
з ядром
$$
\Psi_\beta(t)=\sum\limits_{k=1}^\infty\psi(k)\cos\big(kt-\frac{\beta\pi}{2}\big).
$$
Якщо ж $f\in C_\beta^\psi\mathfrak{N}$, то рівність $(\ref{0.3})$ виконується для всіх $x\in \mathbb{R}$.

У випадку, коли
\begin{equation}\label{S4.0.41}
\psi(k)=e^{-\alpha k^r}, \alpha >0,\ r>0
\end{equation}
множини $L^\psi_\beta\mathfrak{N}$ і $C^\psi _\beta\mathfrak{N}$ будемо позначати через $L^{\alpha ,r}_\beta\mathfrak{N}$ і $C^{\alpha ,r}_\beta\mathfrak{N}$ відповідно. Множини $C^{\alpha ,r}_\beta\mathfrak{N}$ складаються з узагальнених інтегралів Пуассона (див., наприклад, [\ref{Falaleev 2001},~c.~926]), тобто функцій вигляду
\begin{equation}\label{S4.0.42}
f(x)=\frac{a_0}{2}+\frac{1}{\pi }\int\limits_{-\pi }^{\pi }\varphi (x-t)P_{\alpha,r,\beta}(t)\,dt,\ \ \varphi\in \mathfrak{N}, \ \ \varphi\perp1, \ \ a_0\in\mathbb{R},
\end{equation}
де
$$
P_{\alpha,r,\beta}(t)=\sum\limits_{k=1}^{\infty}e^{-\alpha k^r}\cos\big(kt-\frac{\beta \pi }{2}\big)
$$
--- узагальнені ядра Пуассона з параметрами $\alpha >0, \ \ r>0 \ \ \beta \in\mathbb{R}$. Функцію $\varphi$ в рівності (\ref{S4.0.42}) будемо позначати через $f_\beta^{\alpha,r}$. Узагальнені інтеграли Пуассона (тобто функції $f$ вигляду (\ref{S4.0.42})) будемо позначати через ${\cal J}_\beta^{\alpha,r}(\varphi)$.

При $r=1$ ядра $P_{\alpha,r,\beta}(t)$ є звичайними ядрами Пуассона з параметрами $\alpha $ і $\beta$, тобто
$$
P_{\alpha,1,\beta}(t)=\sum\limits_{k=1}^{\infty}e^{-\alpha k}\cos\big(kt-\frac{\beta \pi }{2}\big),
$$
а інтеграли ${\cal J}_\beta^{\alpha,1}(\varphi)$ --- звичайними інтегралами Пуассона функції ${\varphi\in\mathfrak{N}}$.
Неважко переконатися, що при $r=1$ послідовність $\psi(k)$ вигляду (\ref{S4.0.41}) задовольняє умову ${\cal D}_q$, де  $q=e^{-\alpha}$, а при $r>1$ --- умову ${\cal D}_0$.

Одиничну кулю в просторі $L_s, \ 1\leq s\leq\infty,$ позначимо через $U_s$: ${U_s=\big \{ f:f\in L_s, \|f\|_s\leq 1\big\}}$. Для скорочення запису покладемо $C_\beta^\psi U_s=C_{\beta,s}^\psi$, \ $L_\beta^\psi U_s=L_{\beta,s}^\psi$.

Підпростір тригонометричних поліномів $t_{m-1}$, порядок яких не перевищує $m-1$, позначимо через $\mathcal{T}_{2m-1}$. Величина
$$
E_m(f)_{X}=\mathop {\rm \inf}\limits_{t_{m-1}\in{\mathcal{T}}_{2m-1}} \|f-t_{m-1}\|_X
$$
є найкращим наближенням функції $f\in X\subset L_1$ в метриці простору $X$ тригонометричними поліномами порядку $m-1$. Далі в ролі $X$ виступатимуть простори $C$ або $L_s, \  1\leq s\leq\infty$.

Якщо $f\in L_1$, то поліноми вигляду
$$
V_{n,p}(f)=V_{n,p}(f;x)=\frac{1}{p}\sum_{k=n-p}^{n-1}S_k(f;x),
$$
де $S_k(f)=S_k(f;x)$ --- частинні суми Фур'є  порядку $k$ функції $f$, а $p=p(n)$ --- певний натуральний параметр, $p\leqslant n$ називають сумами  Валле  Пуссена функції~$f$.  При $p=1$ суми Валле Пуссена $V_{n,1}(f)$ є  частинними сумами Фур'є $S_{n-1}(f)$ порядку $n-1$, якщо ж $p=n$, то суми
$V_{n,p}(f)$ перетворюються у відомі суми Фейєра $\sigma_{n-1}(f)$
порядку $n-1$:
$$
\sigma_{n-1}(f)=\sigma_{n-1}(f;x)=\frac{1}{n}\sum_{k=0}^{n-1}S_k(f;x).
$$

Для сум Валле Пуссена має місце нерівність Лебега (див., наприклад, [\ref{Stech 1978},~c.~61]):
\begin{equation}\label{neriv}
\|f-V_{n,p}(f)\|_C\leq(\|V_{n,p}\|+1)E_{n-p+1}(f)_C,
\end{equation}
де $\|V_{n,p}\|\mathop{=}\limits^{\rm df}\mathop{\sup}\limits_{\|f\|_C\leq1}\|V_{n,p}(f;\cdot)\|_C$ --- норма оператора $V_{n,p}$. При $p=1$
формула (\ref{neriv}) належить А. Лебегу [\ref{Lebeg1910}]. Оцінки величин $\|V_{n,p}\|$
встановлювались в роботах Валле~Пуссена [\ref{Valle Pyssen 1919}], С.М.~Нікольського~[\,\ref{Nikol 1940}\,], С.Б.~Стєчкіна [\ref{Stech 1951}] та ін. Подальший розвиток вказаної тематики пов'язаний з
дослідженнями О.Д. Габісонії~[\ref{Gabis 1965}], А.А.~Захарова [\ref{Zaharov 1968}], С.Б.~Стєчкіна [\ref{Stech 1978}] та ін. В цих роботах значення величин $\|\rho_{n,p}(f;\cdot)\|_C$, де $\rho_{n,p}(f;x)\mathop{=}\limits^{\rm df}f(x)-V_{n,p}(f;x)$,
оцінювались через найкращі наближення $E_m(f)_C$. Зазначимо, що остаточні порядкові результати в даному напрямі
належать С.Б. Стєчкіну [\ref{Stech 1978}, c. 62], який довів, що для довільної функції $f\in C$ і будь-яких $n, p\in\mathbb{N}$, $p\leq n$ виконується нерівність
\begin{equation}\label{0.7}
\|\rho_{n,p}(f;\cdot)\|_C\leq A
\sum_{k=0}^{n-1}\frac{E_{n-p+k+1}(f)_{C}}{p+k},
\end{equation}
де $A$ --- деяка абсолютна стала. Ним же було доведено, що  нерівність (\ref{0.7}) є точною за порядком для досить широкої множини функціональних компактів. Дійсно, якщо $\varepsilon=\varepsilon_m, \  \  m\in\mathbb{N}$ --- довільна монотонно незростаюча, нескінченно мала послідовність, а $C(\varepsilon)$ --- множина функцій $\varphi\in C$, для яких $E_m(\varphi)_C\leq \varepsilon_m$  $\forall m\in\mathbb{N}$, то  існують абсолютні сталі $A_1$ і $A_2$ такі, що для всіх $n, p\in\mathbb{N}$, $p\leq n$,
\begin{equation}\label{0.8}
A_1\sum_{k=0}^{n-1}\frac{\varepsilon_{n-p+k+1}}{p+k}\leq\mathop{\sup}\limits_{f\in
C(\varepsilon)}\|\rho_{n,p}(f;\cdot)\|_C\leq A_2\sum_{k=0}^{n-1}\frac{\varepsilon_{n-p+k+1}}{p+k}.
\end{equation}
При $p=1$ і $p=n$ співвідношення (\ref{0.8}) доведено в роботах [\ref{Oscolcov 1975}] і [\ref{Stech 1961}] відповідно.

Асимптотична поведінка величин
\begin{equation}\label{st1}
{\cal E}(\mathfrak{N};V_{n,p})_X=\mathop{\sup}\limits_{f\in\mathfrak{N}}\|f-V_{n,p}(f)\|_X, \ \ \ X\subset L_1
\end{equation}
на деяких важливих класах періодичних функцій $\mathfrak{N}\subset X$ при $X=C$ досліджувались багатьма авторами, такими як   А.М.~Колмогоров [\ref{Kolmogorov 1935 }], С.М.~Нікольський~[\ref{Nikol 1941}, \ref{Nikol 1945}], О.П.~Тіман [\ref{Timan 1951 }] та інші. Детальніше ознайомитись з відомими результатами у даному напрямку та з історією питання можна, наприклад, в бібліографічних коментарях до монографій [\ref{Step monog 2002(1)}, \ref{Step,Rycasov monog 2007}, \ref{Step monog 2002(2)}]. На класах $L_\beta^\psi\mathfrak{N}$ та $C_\beta^\psi\mathfrak{N}$ асимптотична поведінка величин (\ref{st1}) при $\psi\in{\cal D}_q$ досліджувалась в роботах [\ref{RykasovHcay2002(1)}--\ref{Serd_Ovsii_2011}]. Для класів $C_{\beta,s}^\psi$ при $\psi\in {\cal D}_q, \ q\in(0,1), \ 1\leq s \leq\infty$, $\beta\in\mathbb{R}$ найбільш остаточні результати по дослідженню асимптотики величин ${\cal E}(C_{\beta,s}^\psi;V_{n,p})$, $n,p\in\mathbb{N}, \ p\leq n$ містяться в роботах [\ref{Serd 2010}, c.~1674] і [\ref{Serd_Ovsii_2011},~с.~4] з яких, зокрема, випливає, що при $n-p\rightarrow\infty$ має місце рівність
\begin{equation}\label{S4.T.1.6(2)}
{\cal E}(C^\psi_{\beta,s};V_{n,p})_C=\frac{\psi(n-p+1)}{p}\bigg(\frac{\|\cos t\|_{s'}}{\pi^{1+1/s'}}K_{q,p}(s')+$$$$+
O(1)\Big(\frac{q\delta(s)}{(n-p+1)(1-q)^{\sigma(s',p)}}+\frac{\varepsilon_{n-p+1}}{(1-q)^2}\min\big\{p,\frac{1}{1-q}\big\}\Big)\bigg),
\end{equation}
де $\frac{1}{s}+\frac{1}{s'}=1$,
\begin{equation}\label{S4.T.1.3}
K_{q,p}(u)=2^{-1/u}\big\|\frac{\sqrt{1-2q^p\cos pt+q^{2p}}}{1-2q\cos t+q^2}\big\|_{u}, \ 1\leq u\leq \infty, q\in(0,1),p\in\mathbb{N},
\end{equation}

\begin{equation}\label{S4.T.1.4}
\sigma(u,p)=\left\{\begin{array}{ll}
1 & {\mbox{\rm при}}\ \ \ \ \ \ \ \ u=1 \ \ \ \ \ {\mbox{\rm і}} \ \ \ \ p=1, \\
2 & {\mbox{\rm при}}\ \ \  1< u\leq\infty \ \ \  {\mbox{\rm і}}  \ \ \ \ p=1, \\
3 & {\mbox{\rm при}}\ \ \ 1\leq u\leq\infty \ \ \ {\mbox{\rm і}} \ \ \ \ p\in\mathbb{N}\backslash\{1\},
\end{array}\right.
\end{equation}

\begin{equation}\label{S4.T.1.2}
\delta(s)= \left\{\begin{array}{ll}
0,& \ \ s=2,\\
1,&  \ \  s\in[1,\infty]\setminus \{2\},
\end{array}\right.
\end{equation}

\begin{equation}\label{S4.T.1.5}
\varepsilon_m=\varepsilon_m(\psi)=\mathop{\sup}\limits_{k\geq m}\left|\frac{\psi(k+1)}{\psi(k)}-q\right|,
\end{equation}
а $O(1)$ --- величина, рівномірно обмежена відносно всіх розглядуваних параметрів.

В роботі [\ref{Stepan i Serduk 2000}] встановлено аналоги нерівностей Лебега для функцій $f$ з ${L_\beta^{\alpha,1}L_s, \ 1\leq s\leq \infty}$, $ \alpha>0$, в яких оцінки відхилень ${\|f(\cdot)-S_{n-1}(f;\cdot)\|_s}$ виражалися через найкращі наближення функцій $f_\beta^{\alpha,1}$  в метриках просторів $L_s$.  А саме, доведено, що для довільних $f\in L_\beta^{\alpha,1}L_s, \alpha>0, \ \beta\in\mathbb{R}, \ 1\leq s\leq\infty$ має місце нерівність
\begin{equation}\label{xxxx}
\|f-S_{n-1}(f)\|_s\leq\Big(\frac{8q^n}{\pi^2}\mathbf{K}(q)+O(1)\frac{q^n}{(1-q)^2n}\Big)E_n(f_\beta^{\alpha,1})_{L_s},
\end{equation}
де $q=e^{-\alpha}$,
$$
\textbf{K}(q)=\int^{\frac{\pi}{2}}_0 \frac{dv}{\sqrt{1-q^2\sin^2v}}
$$
--- повний еліптичний інтеграл першого роду, а  $O(1)$ --- величина, рівномірно обмежена по параметрах $n,q,\beta$, $s$,
$f\in L_\beta^{\alpha,1}L_s $. При $s=1,\infty$ доведено асимптотичну непокращуваність одержаних нерівностей.
В [\ref{Serduk Mysienko 2010}] результати роботи [\ref{Stepan i Serduk 2000}] були узагальнені на випадок сум Валле Пуссена.

Дана робота, яка є продовженням досліджень [\ref{Stepan i Serduk 2000}, \ref{Serduk Mysienko 2010}], присвячена встановленню  асимптотично непокращуваних аналогів нерівностей типу Лебега для сум $V_{n,p}(f)$ на множинах $C_\beta^\psi L_s$ при $\psi\in {\cal D}_q,\,\,\, q\in[0,1)$, $\beta\in\mathbb{R}$ і ${1\leq s\leq\infty}$.

\vskip 0.3cm

\bf Теорема 1. \it{Нехай $\psi\in{\cal D}_q$, $0<q<1$, $\beta\in \mathbb{R}$, $n,p\in\mathbb{N}$, $p\leq n$ і $1\leq s\leq \infty$. Тоді для довільної функції $f\in C_\beta^\psi L_s$ справедлива нерівність
\begin{equation}\label{S4.T.1.1}
\|\rho_{n,p}(f;x)\|_C\leq\frac{\psi(n-p+1)}{p}\bigg(\frac{\|\cos t\|_{s'}}{\pi^{1+1/s'}}K_{q,p}(s')+$$$$+
O(1)\Big(\frac{q\delta(s)}{(n-p+1)(1-q)^{\sigma(s',p)}}+\frac{\varepsilon_{n-p+1}}{(1-q)^2}\min\big\{p,\frac{1}{1-q}\big\}\Big)\bigg)
E_{n-p+1}(f_\beta^\psi)_{L_s},
\end{equation}
в якій $s'=\frac{s}{s-1}$, а $K_{q,p}(s')$, $\sigma(s',p)$, $\delta(s)$ і $\varepsilon_{n-p+1}$ визначаються формулами {\rm(\ref{S4.T.1.3})}, {\rm(\ref{S4.T.1.4})}, {\rm(\ref{S4.T.1.2})} і {\rm(\ref{S4.T.1.5})} відповідно.

При цьому для будь-якої функції $f\in C^\psi_\beta L_s,  1\leq s\leq\infty$ і довільних $n,p\in \mathbb{N}$, $p\leq n$, в множині
$C^\psi_\beta L_s,  1\leq s\leq\infty$, знайдеться функція $F(x)=F(f;n;p;x)$ така, що
$E_{n-p+1}{(F_\beta^\psi)}_{L_s}=E_{n-p+1}{(f_\beta^\psi)}_{L_s}$, і для неї при $n-p\rightarrow\infty$ виконується рівність
\begin{equation}\label{S4.T.1.6}
\|\rho_{n,p}(F;x)\|_C=\frac{\psi(n-p+1)}{p}\bigg(\frac{\|\cos t\|_{s'}}{\pi^{1+1/s'}}K_{q,p}(s')+$$$$+
O(1)\Big(\frac{q\delta(s)}{(n-p+1)(1-q)^{\sigma(s',p)}}+\frac{\varepsilon_{n-p+1}}{(1-q)^2}\min\big\{p,\frac{1}{1-q}\big\}\Big)\bigg)
E_{n-p+1}(F_\beta^\psi)_{L_s}.
\end{equation}
У {\rm(\ref{S4.T.1.1})} і {\rm(\ref{S4.T.1.6})} $O(1)$ --- величини, рівномірно обмежені відносно всіх розглядуваних параметрів.

\rm Із теореми 1 випливає, що нерівність (\ref{S4.T.1.1}) є асимптотично точною при ${n-p\rightarrow\infty}$ на всіх множинах $C_\beta^\psi L_s$ при довільних ${q\in(0,1)}$, $\beta\in\mathbb{R}$ і ${1\leq s\leq\infty}$. Покажемо, що ця ж нерівність залишається асимптотично точною і на деяких важливих підмножинах з $C_\beta^\psi L_s$. Дійсно, розглядаючи точні верхні межі обох частин (\ref{S4.T.1.1}) по класах $C^\psi_{\beta,s}, \ 1\leq s\leq\infty$ і враховуючи, що $E_{n-p+1}(f^\psi_\beta)_{L_s}\leq1$,  отримуємо
\begin{equation}\label{S4.T.1.6(1)}
{\cal E}(C^\psi_{\beta,s};V_{n,p})_C\leq\frac{\psi(n-p+1)}{p}\bigg(\frac{\|\cos t\|_{s'}}{\pi^{1+1/s'}}K_{q,p}(s')+$$$$+
O(1)\Big(\frac{q\delta(s)}{(n-p+1)(1-q)^{\sigma(s',p)}}+\frac{\varepsilon_{n-p+1}}{(1-q)^2}\min\big\{p,\frac{1}{1-q}\big\}\Big)\bigg).
\end{equation}
Співставляючи останнє співвідношення з  асимптотичною рівністю (\ref{S4.T.1.6(2)}), приходимо до висновку, що в (\ref{S4.T.1.6(1)}) можна поставити знак "дорівнює".

Важливими прикладами ядер $\Psi_\beta$ коефіцієнти $\psi(k)$ яких задовольняють умову ${\psi\in {\cal D}_q}$ є

--- полігармонічні ядра Пуассона (див. [\ref{Timan 2009}, c. 256--257])
\begin{equation}\label{S4.T.1.6(9)}
P_{q,\beta}(l,t)=\sum\limits_{k=1}^\infty\psi_l(k)\cos\Big(kt-\frac{\beta\pi}{2}\Big),\,\, l\in\mathbb{N},\,\, \beta\in\mathbb{R},
\end{equation}
де
\begin{equation}\label{S4.T.1.6(10)}
\psi_l(k)=q^k\Big(1+\sum\limits_{j=1}^{l-1}\frac{(1-q^2)^j}{j!2^j}\prod\limits_{\nu=0}^{j-1}(k+2\nu)\Big),\,\, q\in(0,1),
\end{equation}
--- ядра Пуассона для рівняння теплопровідності (див. [\ref{Axiezer monog}, c. 269])
\begin{equation}\label{S4.T.1.6(9')}
{\cal P}_{q,\beta}(t)=\sum\limits_{k=1}^\infty\frac{2}{q^k+q^{-k}}\cos\Big(kt-\frac{\beta\pi}{2}\Big),\,\, \ q\in(0,1), \,\, \beta\in\mathbb{R},
\end{equation}
--- ядра Неймана (див. [\ref{Step monog 2002(1)}, c. 361])
\begin{equation}\label{S4.T.1.6(11)}
N_{q,\beta}(t)=\sum\limits_{k=1}^\infty\frac{q^k}{k}\cos\Big(kt-\frac{\beta\pi}{2}\Big),\,\, q\in(0,1),\,\, \beta\in\mathbb{R}
\end{equation}
та інші.

Для величин $\varepsilon_m$ вигляду (\ref{S4.T.1.5}), що породжуються послідовностями ${\psi(k)=\psi_l(k)}$ ядер $P_{q,\beta}(l,t)$ при $l=1$  справджується  тотожність
\begin{equation}\label{S4.T.1.6(13)}
\varepsilon_m\equiv 0,
\end{equation}
а при $l\in\mathbb{N}\setminus\{1\}$, як доведено в [\ref{Serduk Cay 2011}, c. 108] (див. також [\ref{Serduk Sokolenko 2010}, c. 180]), має місце нерівність
\begin{equation}\label{S4.T.1.6(14)}
\varepsilon_m(l)\leq\frac{(2l-3)q}{m},\,\,\ \ \ m\in\mathbb{N};
\end{equation}
для  $\varepsilon_m$, породжених послідовностями $\psi(k)=\frac{2}{q^k+q^{-k}}$ ядер ${\cal P}_{q,\beta}(t)$, справедлива оцінка
\begin{equation}\label{S4.T.1.6(12')}
\varepsilon_m=\mathop{\sup}\limits_{k\geq m}\frac{q^{2k+1}(1-q^2)}{1+q^{2(k+1)}}=\frac{q^{2m+1}(1-q^2)}{1+q^{2(m+1)}}<q^{2m+1};
\end{equation}
для $\varepsilon_m$, що породжуються послідовностями $\psi(k)=\frac{q^k}{k}$  ядер $N_{q,\beta}(t)$, легко довести рівність
\begin{equation}\label{S4.T.1.6(12)}
\varepsilon_m=\mathop{\sup}\limits_{k\geq m}\frac{q}{k+1}=\frac{q}{m+1}.
\end{equation}

Із теореми 1 і формул (\ref{S4.T.1.6(13)})--(\ref{S4.T.1.6(12)}) отримуємо наступні твердження.

{\textbf{\textit{Наслідок {\rm\bf1}.}}}  \it Нехай множини $C^\psi_\beta L_s,  \ \beta\in \mathbb{R}, \  1\leq s\leq\infty$ породжуються коефіцієнтами $\psi(k)=\psi_l(k)$ ядер $P_{q,\beta}(l,t)$, $l\in\mathbb{N}$, вигляду {\rm(\ref{S4.T.1.6(9)})}. Тоді для довільних  $f\in C_\beta^\psi L_s$, $n,p\in\mathbb{N}$, $p\leq n$ справедлива нерівність
\begin{equation}\label{S4.T.1.6(15)}
\|\rho_{n,p}(f;x)\|_C\leq\frac{q^{n-p+1}}{p}\big(1+\sum\limits_{j=1}^{l-1}\frac{(1-q^2)^j}{j!2^j}\prod\limits_{\nu=0}^{j-1}(k+2\nu)\big)\times$$$$\times
\bigg(\frac{\|\cos t\|_{s'}}{\pi^{1+1/s'}}K_{q,p}(s')+
O(1)\Big(\frac{q\delta(s)}{(n-p+1)(1-q)^{\sigma(s',p)}}+$$$$+\frac{lq}{(n-p+1)(1-q)^2}\min\big\{p,\frac{1}{1-q}\big\}\Big)\bigg)
E_{n-p+1}(f_\beta^\psi)_{L_s},
\end{equation}
де
$s'=\frac{s}{s-1}$, а $K_{q,p}(s')$, $\sigma(s',p)$ і $\delta(s)$ визначаються формулами {\rm(\ref{S4.T.1.3})}, {\rm(\ref{S4.T.1.4})} і {\rm(\ref{S4.T.1.2})} відповідно.

При цьому для будь-якої функції $f\in C^\psi_\beta L_s,  1\leq s\leq\infty$ і довільних $n,p\in \mathbb{N}$, $p\leq n$, в множині
$C^\psi_\beta L_s,  1\leq s\leq\infty$, знайдеться функція $F(x)=F(f;n;p;x)$ така, що $E_{n-p+1}{(F_\beta^\psi)}_{L_s}=E_{n-p+1}{(f_\beta^\psi)}_{L_s}$, і для неї при $n-p\rightarrow\infty$ виконується рівність
\begin{equation}\label{S4.T.1.6(15')}
\|\rho_{n,p}(F;x)\|_C=\frac{q^{n-p+1}}{p}\big(1+\sum\limits_{j=1}^{l-1}\frac{(1-q^2)^j}{j!2^j}\prod\limits_{\nu=0}^{j-1}(k+2\nu)\big)\times$$$$\times
\bigg(\frac{\|\cos t\|_{s'}}{\pi^{1+1/s'}}K_{q,p}(s')+
O(1)\Big(\frac{q\delta(s)}{(n-p+1)(1-q)^{\sigma(s',p)}}+$$$$+\frac{lq}{(n-p+1)(1-q)^2}\min\big\{p,\frac{1}{1-q}\big\}\Big)\bigg)
E_{n-p+1}(F_\beta^\psi)_{L_s}.
\end{equation}
У {\rm(\ref{S4.T.1.6(15)})} і {\rm(\ref{S4.T.1.6(15')})} $O(1)$ --- величини, рівномірно обмежені відносно всіх розглядуваних параметрів.

{\textbf{\textit{Наслідок {\rm\bf2}.}}}  \it Нехай множини $C^\psi_\beta L_s, \ \beta\in \mathbb{R}, \ 1\leq s\leq\infty$ породжуються коефіцієнтами $\psi(k)=\frac{2}{q^k+q^{-k}}$ ядер ${\cal P}_{q,\beta}(t)$, вигляду {\rm(\ref{S4.T.1.6(9')})}. Тоді для довільних $f\in C_\beta^\psi L_s$, $n,p\in\mathbb{N}$, $p\leq n$ справедлива нерівність
\begin{equation}\label{S4.T.1.6(14')}
\|\rho_{n,p}(f;x)\|_C\leq$$$$\leq\frac{2q^{n-p+1}}{(1+q^{2(n-p+1)})p}\bigg(\frac{\|\cos t\|_{s'}}{\pi^{1+1/s'}}K_{q,p}(s')+
O(1)\Big(\frac{q\delta(s)}{(n-p+1)(1-q)^{\sigma(s',p)}}+$$$$+\frac{q^{2(n-p)+3}}{(1-q)^2}\min\big\{p,\frac{1}{1-q}\big\}\Big)\bigg)
E_{n-p+1}(f_\beta^\psi)_{L_s},
\end{equation}
де
$s'=\frac{s}{s-1}$, а $K_{q,p}(s')$, $\sigma(s',p)$ і $\delta(s)$ визначаються формулами {\rm(\ref{S4.T.1.3})}, {\rm(\ref{S4.T.1.4})} і {\rm(\ref{S4.T.1.2})} відповідно.

При цьому для будь-якої функції $f\in C^\psi_\beta L_s,  1\leq s\leq\infty$ і довільних $n,p\in \mathbb{N}$, $p\leq n$, в множині
$C^\psi_\beta L_s,  1\leq s\leq\infty$, знайдеться функція $F(x)=F(f;n;p;x)$ така, що $E_{n-p+1}{(F_\beta^\psi)}_{L_s}=E_{n-p+1}{(f_\beta^\psi)}_{L_s}$, і для неї при $n-p\rightarrow\infty$ виконується рівність
\begin{equation}\label{S4.T.1.6(14'')}
\|\rho_{n,p}(F;x)\|_C=$$$$=\frac{2q^{n-p+1}}{(1+q^{2(n-p+1)})p}\big(\frac{\|\cos t\|_{s'}}{\pi^{1+1/s'}}K_{q,p}(s')+
O(1)\big(\frac{q\delta(s)}{(n-p+1)(1-q)^{\sigma(s',p)}}+$$$$+\frac{q^{2(n-p)+3}}{(1-q)^2}\min\big\{p,\frac{1}{1-q}\big\}\big)\big)
E_{n-p+1}(F_\beta^\psi)_{L_s}.
\end{equation}
У {\rm(\ref{S4.T.1.6(14')})} і {\rm(\ref{S4.T.1.6(14'')})} $O(1)$ --- величини, рівномірно обмежені відносно всіх розглядуваних параметрів.

{\textbf{\textit{Наслідок {\rm\bf3}.}}}  \it Нехай множини $C^\psi_\beta L_s, \ \beta\in \mathbb{R}, \ 1\leq s\leq\infty$ породжуються коефіцієнтами $\psi(k)=\frac{q^k}{k}$ ядер $N_{q,\beta}(t)$, вигляду {\rm(\ref{S4.T.1.6(11)})}. Тоді для довільних ${f\in C_\beta^\psi L_s}$, $n,p\in\mathbb{N}$, $p\leq n$ справедлива нерівність
\begin{equation}\label{S4.T.1.6(14)}
\|\rho_{n,p}(f;x)\|_C\leq$$$$\leq\frac{q^{n-p+1}}{p(n-p+1)}\bigg(\frac{\|\cos t\|_{s'}}{\pi^{1+1/s'}}K_{q,p}(s')+
O(1)\Big(\frac{q\delta(s)}{(n-p+1)(1-q)^{\sigma(s',p)}}+$$$$+\frac{q}{(n-p+2)(1-q)^2}\min\big\{p,\frac{1}{1-q}\big\}\Big)\bigg)
E_{n-p+1}(f_\beta^\psi)_{L_s},
\end{equation}
де
$s'=\frac{s}{s-1}$, а $K_{q,p}(s')$, $\sigma(s',p)$ і $\delta(s)$ визначаються формулами {\rm(\ref{S4.T.1.3})}, {\rm(\ref{S4.T.1.4})} і {\rm(\ref{S4.T.1.2})} відповідно.

При цьому для будь-якої функції $f\in C^\psi_\beta L_s,  1\leq s\leq\infty$ і довільних $n,p\in \mathbb{N}$, $p\leq n$, в множині
$C^\psi_\beta L_s,  1\leq s\leq\infty$, знайдеться функція $F(x)=F(f;n;p;x)$ така, що $E_{n-p+1}{(F_\beta^\psi)}_{L_s}=E_{n-p+1}{(f_\beta^\psi)}_{L_s}$, і для неї при $n-p\rightarrow\infty$ виконується рівність
\begin{equation}\label{S4.T.1.6(16)}
\|\rho_{n,p}(F;x)\|_C=$$$$=\frac{q^{n-p+1}}{p(n-p+1)}\bigg(\frac{\|\cos t\|_{s'}}{\pi^{1+1/s'}}K_{q,p}(s')+
O(1)\Big(\frac{q\delta(s)}{(n-p+1)(1-q)^{\sigma(s',p)}}+$$$$+\frac{q}{(n-p+2)(1-q)^2}\min\big\{p,\frac{1}{1-q}\big\}\Big)\bigg)
E_{n-p+1}(F_\beta^\psi)_{L_s}.
\end{equation}
У {\rm(\ref{S4.T.1.6(14)})} і {\rm(\ref{S4.T.1.6(16)})} $O(1)$ --- величини, рівномірно обмежені відносно всіх розглядуваних параметрів.

\rm Оскільки, як зазначалося вище, при $p=1$ суми Валле Пуссена $V_{n,p}(f)$ перетворюються в суми Фур'є $S_{n-1}(f)$, то з
(\ref{S4.T.1.1}) для довільної $f\in C_\beta^\psi L_s$, $\psi\in{\cal D}_q$, $q\in(0,1)$, $1\leq s\leq\infty$  випливає нерівність
\begin{equation}\label{S4.T.1.6(3)}
\|f(x)-S_{n-1}(f;x)\|_C\leq$$$$\leq\psi(n)\bigg(\frac{\|\cos t\|_{s'}}{\pi^{1+1/s'}}K_{q,1}(s')+
O(1)\Big(\frac{q\delta(s)}{n(1-q)^{\sigma(s',1)}}+\frac{\varepsilon_{n}}{(1-q)^2}\Big)\bigg)E_{n}(f_\beta^\psi)_{L_s}.
\end{equation}

В роботі [\ref{Serduk 2012}] (див. формулу (25)) було показано, що при $1\leq u<\infty$
$$
K_{q,1}(u)=\pi^{\frac{1}{u}}F^{\frac{1}{u}}\left(\frac{u}{2},\frac{u}{2};1;q^2\right),
$$
де
$$
F(a,b;c;z)=1+\sum\limits_{k=1}^\infty\frac{(a)_k(b)_k}{(c)_k}\frac{z^k}{k!}, \ \ (x)_k=x(x+1)(x+2)\cdots(x+k-1)
$$
--- гіпергеометрична функція Гаусса.
Тому для $1<s\leq\infty$ нерівність (\ref{S4.T.1.6(3)}) можемо переписати у вигляді
\begin{equation}\label{S4.T.1.6(4)}
\|f(x)-S_{n-1}(f;x)\|_C\leq$$$$\leq\psi(n)\bigg(\frac{\|\cos t\|_{s'}}{\pi}F^{\frac{1}{s'}}\big(\frac{s'}{2},\frac{s'}{2};1;q^2\big)+
O(1)\Big(\frac{q\delta(s)}{n(1-q)^{\sigma(s',1)}}+\frac{\varepsilon_{n}}{(1-q)^2}\Big)\bigg)E_{n}(f_\beta^\psi)_{L_s}.
\end{equation}
Якщо $s'=1$, то (див. [\ref{Gren}, c. 919])
\begin{equation}\label{S4.T.1.6(5)}
F\Big(\frac{1}{2},\frac{1}{2};1;q^2\Big)=2{\bf{K}}(q),
\end{equation}
де ${\bf{K}}(q)$ --- повний еліптичний інтеграл першого роду.

Із (\ref{S4.T.1.6(4)}), з урахуванням (\ref{S4.T.1.4}), (\ref{S4.T.1.2})  і (\ref{S4.T.1.6(5)}),  для довільної $f\in C_\beta^\psi L_\infty,  \ \beta\in \mathbb{R}$, $\psi\in{\cal D}_q$, $q\in(0,1)$ можемо записати нерівність
\begin{equation}\label{S4.T.1.6(7)}
\|f(x)-S_{n-1}(f;x)\|_C\leq$$$$\leq\psi(n)\bigg(\frac{8}{\pi^2}{\bf{K}}(q)+
O(1)\Big(\frac{q}{n(1-q)}+\frac{\varepsilon_{n}}{(1-q)^2}\Big)\bigg)E_{n}(f_\beta^\psi)_{L_\infty}.
\end{equation}
Якщо $\psi(k)=e^{-\alpha k},\ \ \alpha>0,$ то з (\ref{S4.T.1.6(7)}) випливає нерівність (\ref{xxxx}) при $s=\infty$.

Розглядаючи точні верхні межі обох частин (\ref{S4.T.1.6(7)}) по класах $C^\psi_{\beta,\infty}$, отримуємо
\begin{equation}\label{S4.T.1.6(4')}
{\cal E}\big(C^\psi_{\beta,\infty};S_{n-1}\big)_C\leq\psi(n)\bigg(\frac{8}{\pi^2}\mathbf{K}(q)+
O(1)\Big(\frac{q}{n(1-q)}+\frac{\varepsilon_{n}}{(1-q)^2}\Big)\bigg),
\end{equation}
Співставляючи це співвідношення з отриманою в роботі [\ref{Serd_Step 2000},~с.~384] асимптотичною при $n\rightarrow\infty$ рівністю
\begin{equation}\label{S4.T.1.6(6)}
{\cal E}\big(C^\psi_{\beta,\infty};S_{n-1}\big)_C=\psi(n)\bigg(\frac{8}{\pi^2}\mathbf{K}(q)+
O(1)\Big(\frac{q}{n(1-q)}+\frac{\varepsilon_{n}}{(1-q)^2}\Big)\bigg),
\end{equation}
приходимо до висновку, що при $s=\infty$ в співвідношенні (\ref{S4.T.1.6(4')}) можна поставити знак "дорівнює".

При $s'=2$, як легко переконатися,
\begin{equation}\label{S4.T.1.6(7')}
F^{\frac{1}{s'}}\Big(\frac{s'}{2},\frac{s'}{2};1;q^2\Big)=F^{\frac{1}{2}}(1,1;1;q^2)=\frac{1}{\sqrt{1-q^2}}.
\end{equation}
Із (\ref{S4.T.1.6(4)}), з урахуванням {\rm(\ref{S4.T.1.4})}, (\ref{S4.T.1.2}) і (\ref{S4.T.1.6(7')}), для довільних $f\in C_\beta^\psi L_2, \ \beta\in \mathbb{R}$, $\psi\in{\cal D}_q$, $q\in(0,1)$ можемо записати нерівність
\begin{equation}\label{S4.T.1.6(8)}
\|f(x)-S_{n-1}(f;x)\|_C\leq\psi(n)\bigg(\frac{1}{\sqrt{\pi(1-q^2)}} +O(1)\frac{\varepsilon_n}{(1-q)^2}\bigg)E_{n}(f_\beta^\psi)_{L_2}.
\end{equation}
Розглядаючи точні верхні межі обох частин нерівності (\ref{S4.T.1.6(8)}) по класах $C^\psi_{\beta,2}$, отримуємо
\begin{equation}\label{S4.T.1.6(8')}
{\cal E}\big(C^\psi_{\beta,2};S_{n-1}\big)_C\leq\psi(n)\bigg(\frac{1}{\sqrt{\pi(1-q^2)}} +O(1)\frac{\varepsilon_n}{(1-q)^2}\bigg).
\end{equation}
Як випливає з формули (30) роботи [\ref{Stepan i Serduk 2000}] і формули (65$'$) роботи [\ref{Serd 2005}] в співвідношенні (\ref{S4.T.1.6(8')}) можна поставити знак "дорівнює". Зазначимо також, що в [\ref{Serduk Sokolenko 2011}] встановлені точні рівності величин ${\cal E}(C_{\beta,2}^\psi;S_{n-1})$ за умови $\sum\limits_{k=1}^\infty\psi^2(k)<\infty$.

При $s'=\infty$, як випливає з (\ref{S4.T.1.3}),
\begin{equation}\label{S4.T.1.6(81)}
K_{q,1}(s')=K_{q,1}(\infty)=\Big\|\frac{1}{\sqrt{1-2q\cos t+q^2}}\Big\|_\infty=\frac{1}{(1-q)}.
\end{equation}
Із (\ref{S4.T.1.6(3)}), з урахуванням {\rm(\ref{S4.T.1.4})}, (\ref{S4.T.1.2}) і (\ref{S4.T.1.6(81)}), для довільної $f\in C_\beta^\psi L_1 , \ \beta\in \mathbb{R}$, $\psi\in{\cal D}_q$, $q\in(0,1)$ можемо записати нерівність
\begin{equation}\label{S4.T.1.6(82)}
\|f(x)-S_{n-1}(f;x)\|_C\leq\psi(n)\bigg(\frac{1}{\pi(1-q)}+O(1)\frac{\varepsilon_n+q/n}{(1-q)^2}\bigg)E_{n}(f_\beta^\psi)_{L_1}.
\end{equation}
Розглядаючи точні верхні межі обох частин нерівності (\ref{S4.T.1.6(82)}) по класах $C^\psi_{\beta,1}$, отримуємо
\begin{equation}\label{S4.T.1.6(83)}
{\cal E}\big(C^\psi_{\beta,1};S_{n-1}\big)_C\leq\psi(n)\bigg(\frac{1}{\pi(1-q)}+O(1)\frac{\varepsilon_n+q/n}{(1-q)^2}\bigg).
\end{equation}
Як випливає з формули (30) роботи [\ref{Stepan i Serduk 2000}] і формули (63) роботи [\ref{Serd 2005}] в (\ref{S4.T.1.6(83)}) при $n\rightarrow\infty$ можна поставити знак "дорівнює".

{\textbf{\textit{Доведення теореми {\rm\bf1}.}}} \ \rm
Нехай $f\in C_\beta^\psi L_s$, $1\leq s\leq\infty$. Тоді в кожній точці $x\in\mathbb{R}$ (див., наприклад, [\ref{Rykasov 2003}, c. 810]) має місце інтегральне зображення
\begin{equation}\label{S4.T.1.7}
\rho_{n,p}(f;x)=f(x)-V_{n,p}(f;x)=\frac{1}{\pi}\int\limits_{-\pi}^\pi f_\beta^\psi(x-t)\Psi_{1,n,p}(t)dt,
\end{equation}
в якому
\begin{equation}\label{S4.T.1.8}
\Psi_{j,n,p}(t)=\sum\limits_{k=n-p+j}^\infty\tau_{n,p}(k)\psi(k)\cos\Big(kt-\frac{\beta\pi}{2}\Big), \ j\in\mathbb{N},
\end{equation}
а $\tau_{n,p}(k)$ визначається наступним чином:
\begin{equation}\label{S4.T.1.9}
\tau_{n,p}(k)= \left\{\begin{array}{ll}
1-\frac{n-k}{p},& \ \ n-p+1\leq k\leq n-1,\\
1,&  \ \  k\geq n.
\end{array}\right.
\end{equation}

Покладемо
\begin{equation}\label{S4.T.1.13}
r_{n,p}(t)=\sum_{k=n-p+2}^\infty\tau_{n,p}(k)\Big(\frac{\psi(k)}{\psi(n-p+1)}-\frac{q^k}{q^{n-p+1}}\Big)
\cos\Big(kt-\frac{\beta\pi}{2}\Big).
\end{equation}
Враховуючи (\ref{S4.T.1.13}), перепишемо (\ref{S4.T.1.7}) у вигляді
$$
\rho_{n,p}(f;x)=$$$$=\psi(n-p+1)\bigg(\frac{q^{-(n-p+1)}}{\pi}\int\limits_{-\pi}^\pi f_\beta^\psi(x-t)\sum\limits_{k=n-p+1}^\infty\tau_{n,p}(k)q^k\cos\Big(kt-\frac{\beta\pi}{2}\Big)dt+$$$$+
\frac{1}{\pi}\int\limits_{-\pi}^\pi f_\beta^\psi(x-t)r_{n,p}(t)dt\bigg).
$$
Провівши елементарні перетворення, з урахуванням (\ref{S4.T.1.9}), одержимо
$$
\sum\limits_{k=n-p+1}^\infty\tau_{n,p}(k)q^k\cos\Big(kt-\frac{\beta\pi}{2}\Big)=\frac{1}{p}\sum\limits_{k=n-p}^{n-1}\sum\limits_{j=k+1}^\infty q^j\cos\Big(jt-\frac{\beta\pi}{2}\Big).
$$
В роботі [\ref{Serd 2004}, c. 99--100] було показано, що
$$
\sum\limits_{k=n-p}^{n-1}\sum\limits_{j=k+1}^\infty q^j\cos\big(jt-\frac{\beta\pi}{2}\big)=Z_q(t) P_{q,\beta,n,p}(t),
$$
де
$$
Z_q(t)=\frac{1}{\sqrt{1-2q\cos{t}+q^2 }},
$$
$$
P_{q,\beta,n,p}(t)=\sum_{k=n-p+1}^nq^k\cos\Big(kt+\theta_q(t)-\frac{\beta\pi}{2}\Big),
$$
$$
\theta_q (t)=\arctg{\frac{q\sin t}{1-q\cos{t}}}.
$$

Функції $Z_q(t) P_{q,\beta,n,p}(t)$ і $r_{n,p}(t)$ ортогональні до будь-якого тригонометричного полінома $t_{n-p}$ порядок, якого не перевищує $n-p$. Тому для довільного полінома $t_{n-p}$ з ${\cal T}_{2(n-p)+1}$
\begin{equation}\label{S4.T.1.15}
\rho_{n,p}(f;x)=\frac{\psi(n-p+1)}{\pi}\bigg(\int\limits_{-\pi}^\pi\delta_{n,p}(x-t)\Big(\frac{q^{-(n-p+1)}}{p}Z_q(t)P_{q,\beta,n,p}(t)+r_{n,p}(t)\Big)dt\bigg),
\end{equation}
де
\begin{equation}\label{S4.T.1.16}
\delta_{n,p}(\cdot)=f_\beta^\psi(\cdot)-t_{n-p}(\cdot).
\end{equation}
В [\ref{Serd_Ovsii_2011}, c. 6-7] встановлено, що при $n-p\rightarrow\infty$ має місце рівномірна по $t,n,p,q,\psi$ і $\beta$ оцінка
\begin{equation}\label{S4.T.1.17}
|r_{n,p}(t)|=O(1)\frac{\varepsilon_{n-p+1}}{(1-q)^2}\min\big\{1,\frac{1}{p(1-q)}\big\}.
\end{equation}
Враховуючи (\ref{S4.T.1.17}), із рівності (\ref{S4.T.1.15}) одержуємо
\begin{equation}\label{S4.T.1.18}
\rho_{n,p}(f;x)=\frac{\psi(n-p+1)}{\pi p}\int\limits_{-\pi}^\pi\delta_{n,p}(x-t)\Big(q^{-(n-p+1)}Z_q(t)P_{q,\beta,n,p}(t)+$$$$+
O(1)\frac{\varepsilon_{n-p+1}}{(1-q)^2}\min\big\{p,\frac{1}{1-q}\big\}\Big)dt.
\end{equation}
Далі, обравши в (\ref{S4.T.1.18}) в ролі $t_{n-p}(\cdot)$ поліном $t_{n-p}^\ast(\cdot)$ найкращого наближення в просторі $L_s$ функції $f_\beta^\psi(\cdot)$, тобто такий, що
$$
\|f_\beta^\psi-t_{n-p}^\ast\|_s=E_{n-p+1}(f_\beta^\psi)_{L_s}, \ \  1\leq s\leq\infty
$$
і застосувавши нерівність
\begin{equation}\label{S4.T.1.18'}
\big\|\int\limits_{-\pi}^\pi K(t-u)\varphi(u)du\big\|_C\leq \|K\|_{s'}\|\varphi\|_s, \ \varphi\in L_s, \ K\in L_{s'}, \ 1\leq s\leq\infty, \ \frac{1}{s}+\frac{1}{s'}=1,
\end{equation}
(див., наприклад, [\ref{Korn}, c. 43]), маємо
\begin{equation}\label{S4.T.1.19}
\|\rho_{n,p}(f;x)\|_C\leq\frac{\psi(n-p+1)}{\pi p}\Big(q^{-(n-p+1)}\|Z_q(t)P_{q,\beta,n,p}(t)\|_{s'}+$$$$+
O(1)\frac{\varepsilon_{n-p+1}}{(1-q)^2}\min\big\{p,\frac{1}{1-q}\big\}\Big)E_{n-p+1}(f_\beta^\psi)_{L_s}.
\end{equation}

В силу співвідношення (69) роботи [\ref{Serd 2010}] для  $s', \ 1\leq s'\leq\infty$ при $n\rightarrow\infty$ має місце асимптотична рівність
\begin{equation}\label{S4.T.1.20}
 \|Z_q(t)P_{q,\beta,n,p}(t)\|_{s'}=$$$$=q^{n-p+1}\bigg(\frac{\|\cos
t\|_{s'}}{\pi^{1/s'}}K_{q,p}(s')+O(1)\frac{q}{(n-p+1)(1-q)^{\sigma(s',p)}}\bigg).
 \end{equation}
При $s=2$ (див. [\ref{Serd 2010}, c. 1680]) рівність (\ref{S4.T.1.20}) можна покращити оскільки
\begin{equation}\label{S4.T.1.20'}
\|Z_q(t)P_{q,\beta,n,p}(t)\|_{2}=\big\|\sum\limits_{k=n-p}^{n-1}\sum\limits_{j=k+1}^\infty q^j\cos\big(jt-\frac{\beta\pi}{2}\big)\big\|_2=$$$$=\sqrt{\pi}q^{n-p+1}\sqrt{\frac{1+q^2-q^{2p}(2p+1-q^2(2p-1))}{(1-q^2)^3}}=$$$$=
q^{n-p+1}K_{p,q}(2)=q^{n-p+1}\frac{\|\cos t\|_2}{\sqrt{\pi}}K_{p,q}(2).
\end{equation}

Поєднавши нерівність (\ref{S4.T.1.19}) з оцінками (\ref{S4.T.1.20}) і (\ref{S4.T.1.20'}), одержуємо (\ref{S4.T.1.1}).

Доведемо другу частину теореми. Для цього досить показати, що для довільної функції $\varphi\in L_s$, можна вказати функцію $\Phi(\cdot)=\Phi(\varphi;\cdot)$, для якої $E_{n-p+1}(\Phi)_{L_s}=E_{n-p+1}(\varphi)_{L_s}$ при всіх $n,p\in\mathbb{N}, \ p\leq n$, і крім того, має місце рівність
\begin{equation}\label{48}
|\rho_{n,p}(\Phi;0)|=\frac{\psi(n-p+1)}{p}\big(\frac{\|\cos t\|_{s'}}{\pi^{1+1/s'}}K_{q,p}(s')+$$$$+
O(1)\big(\frac{q\delta(s)}{(n-p+1)(1-q)^{\sigma(s',p)}}+\frac{\varepsilon_{n-p+1}}{(1-q)^2}\min\big\{p,\frac{1}{1-q}\big\}\big)\big)
E_{n-p+1}(\Phi_\beta^\psi)_{L_s}.
\end{equation}
В силу, інтегрального зображення (\ref{S4.T.1.15}), оцінки (\ref{S4.T.1.17}) та факту ортогональності функції $Z_q(t)P_{q,\beta,n,p}(t)$ до будь-якого тригонометричного полінома ${t_{n-p}\in{\cal T}_{2(n-p)+1}}$, для довільної функції $f$ з множини $C_\beta^\psi L_s$, $\psi\in {\cal D}_q$, $0<q<1$, виконується рівність
\begin{equation}\label{S4.T.1.21}
|\rho_{n,p}(f;x)|=\psi(n-p+1)\Big(q^{-(n-p+1)}\big|\frac{1}{\pi p}\int\limits_{-\pi}^\pi f_\beta^\psi(x-t) Z_q(t)P_{q,\beta,n,p}(t)dt\big|+$$$$+O(1)\frac{\varepsilon_{n-p+1}}{(1-q)^2}\min\big\{1,\frac{1}{p(1-q)}\big\}E_{n-p+1}(f_\beta^\psi)_{L_s}\Big)=$$$$=
\psi(n-p+1)\Big(q^{-(n-p+1)}|\rho_{n,p}({\cal J}_\beta^{\alpha,1}(f_\beta^\psi);x)|+$$$$+O(1)\frac{\varepsilon_{n-p+1}}{(1-q)^2}\min\big\{1,\frac{1}{p(1-q)}\big\}E_{n-p+1}(f_\beta^\psi)_{L_s}\Big),
\end{equation}
де $\alpha=\ln\frac{1}{q}$.

Внаслідок теореми 3 роботи [\ref{Serduk Mysienko 2010}, c. 313] для довільної функції ${\varphi\in L_s, \ 1\leq s\leq\infty}$ знайдеться функція $\overline{\varphi}(t)$ така, що
\begin{equation}\label{S4.T.1.23}
E_{n-p+1}(\overline{\varphi})_{L_s}=E_{n-p+1}(\varphi)_{L_s}
\end{equation}
і для неї виконується рівність
\begin{equation}\label{S4.T.1.23'}
|\rho_{n,p}({\cal J}_\beta^{\alpha,1}(\overline{\varphi});0)|=$$$$=\frac{q^{n-p+1}}{p}\Big(\frac{\|\cos t\|_{s'}}{\pi^{1+1/s'}}K_{q,p}(s')+
O(1)\frac{q\delta(s)}{(n-p+1)(1-q)^{\sigma(s',p)}}\Big)E_{n-p+1}(\overline{\varphi})_{L_s}.
\end{equation}

Функція $F={\cal J}_\beta^\psi(\overline{\varphi})$ є шуканою, оскільки для неї в силу (\ref{S4.T.1.21})--(\ref{S4.T.1.23'}) виконується (\ref{48}). Теорему 1 доведено.

\bf Теорема 2. \it{Нехай $\psi\in {\cal D}_q$, $0<q<1$, $\beta\in \mathbb{R}$, $n,p\in\mathbb{N}$, $p\leq n$ і $1\leq s\leq \infty$. Тоді для довільної функції $f\in L_\beta^\psi L_1$ справедлива нерівність
\begin{equation}\label{S4.T.1'.1}
\|\rho_{n,p}(f;x)\|_{L_s}\leq\frac{\psi(n-p+1)}{p}\bigg(\frac{\|\cos t\|_{s}}{\pi^{1+1/s}}K_{q,p}(s)+$$$$+
O(1)\Big(\frac{q}{(n-p+1)(1-q)^{\sigma(s,p)}}+\frac{\varepsilon_{n-p+1}}{(1-q)^2}\min\big\{p,\frac{1}{1-q}\big\}\Big)\bigg)
E_{n-p+1}(f_\beta^\psi)_{L_1},
\end{equation}
де $K_{q,p}(s)$ і $\sigma(s,p)$ визначаються формулами {\rm(\ref{S4.T.1.3})} і {\rm(\ref{S4.T.1.4})} відповідно, а $O(1)$ --- величина, рівномірно обмежена відносно всіх розглядуваних параметрів.

{\textbf{\textit{Доведення теореми {\rm\bf2}.}}} \ \rm Нехай $f\in L_\beta^\psi L_1$. Внаслідок інтегрального зображення (\ref{S4.T.1.18}) та твердження 1.5.5 із роботи [\ref{Korn}, с. 43], отримуємо
\begin{equation}\label{S4.T.1'.2}
\|\rho_{n,p}(f;x)\|_{L_s}\leq\frac{\psi(n-p+1)}{\pi p}\Big(q^{-(n-p+1)}\|Z_q(t)P_{q,\beta,n,p}(t)\|_s+$$$$+
O(1)\frac{\varepsilon_{n-p+1}}{(1-q)^2}\min\big\{p,\frac{1}{1-q}\big\}\Big)\|\delta_{n,p}(x-t)\|_1, \ \ 1\leq s\leq\infty,
\end{equation}
де $\delta_{n,p}$ визначається рівністю (\ref{S4.T.1.16}). Обравши в (\ref{S4.T.1'.2}) в якості $t_{n-p}$ поліном  найкращого наближення $t_{n-p}^\ast$ функції $f$ в метриці простору $L_1$, одержимо
\begin{equation}\label{S4.T.1'.3}
\|\rho_{n,p}(f;x)\|_{L_s}\leq\frac{\psi(n-p+1)}{\pi p}\Big(q^{-(n-p+1)}\|Z_q(t)P_{q,\beta,n,p}(t)\|_s+$$$$+
O(1)\frac{\varepsilon_{n-p+1}}{(1-q)^2}\min\big\{p,\frac{1}{1-q}\big\}\Big)E_{n-p+1}(f_\beta^\psi)_{L_1}.
\end{equation}
В силу рівності (69) роботи [\ref{Serd 2010}] маємо
\begin{equation}\label{S4.T.1'.4}
\|Z_q(t)P_{q,\beta,n,p}(t)\|_{s}=$$$$=q^{n-p+1}\bigg(\frac{\|\cos t\|_{s}}{\pi^{1/s}}K_{q,p}(s)+O(1)\frac{q}{(n-p+1)(1-q)^{\sigma(s,p)}}\bigg).
 \end{equation}
З (\ref{S4.T.1'.3}) і (\ref{S4.T.1'.4}) випливає (\ref{S4.T.1'.1}). Теорему 2 доведено.

З теореми 2 і оцінок (\ref{S4.T.1.6(13)})--(\ref{S4.T.1.6(12)}) отримуємо наступні твердження.

{\textbf{\textit{Наслідок {\rm\bf4}.}}}   \it Нехай множини $L^\psi_\beta L_1, \ \beta\in \mathbb{R}$ породжуються коефіцієнтами $\psi(k)=\psi_l(k)$ ядер $P_{q,\beta}(l,t)$, $l\in\mathbb{N}$ вигляду {\rm(\ref{S4.T.1.6(9)})}. Тоді для довільних $f\in L^\psi_\beta L_1$, $n,p\in\mathbb{N}$, $p\leq n, \ 1\leq s\leq\infty$ справедлива нерівність
$$
\|\rho_{n,p}(f;x)\|_{L_s}\leq\frac{q^{n-p+1}}{p}\big(1+\sum\limits_{j=1}^{l-1}\frac{(1-q^2)^j}{j!2^j}\prod\limits_{\nu=0}^{j-1}(k+2\nu)\big)\times$$$$\times
\bigg(\frac{\|\cos t\|_{s}}{\pi^{1+1/s}}K_{q,p}(s)+
O(1)\Big(\frac{q}{(n-p+1)(1-q)^{\sigma(s,p)}}+$$$$+\frac{lq}{(n-p+1)(1-q)^2}\min\big\{p,\frac{1}{1-q}\big\}\Big)\bigg)
E_{n-p+1}(f_\beta^\psi)_{L_1},
$$
де
$K_{q,p}(s)$ і $\sigma(s,p)$  визначаються формулами {\rm(\ref{S4.T.1.3})} і {\rm(\ref{S4.T.1.4})} відповідно, а $O(1)$ --- величина, рівномірно обмежена відносно всіх розглядуваних параметрів.

{\textbf{\textit{Наслідок {\rm\bf5}.}}}  \it Нехай множини $L^\psi_\beta L_1, \ \beta\in \mathbb{R}$ породжуються коефіцієнтами $\psi(k)=\frac{2}{q^k+q^{-k}}$ ядер ${\cal P}_{q,\beta}(t)$ вигляду {\rm(\ref{S4.T.1.6(9')})}. Тоді для довільних $f\in L^\psi_\beta L_1$, $n,p\in\mathbb{N}$, $p\leq n, \ 1\leq s\leq\infty $ справедлива нерівність
$$
\|\rho_{n,p}(f;x)\|_{L_s}\leq$$$$\leq\frac{2q^{n-p+1}}{(1+q^{2(n-p+1)})p}\bigg(\frac{\|\cos t\|_{s}}{\pi^{1+1/s}}K_{q,p}(s)+
O(1)\Big(\frac{q}{(n-p+1)(1-q)^{\sigma(s,p)}}+$$$$+\frac{q^{2(n-p)+3}}{(1-q)^2}\min\big\{p,\frac{1}{1-q}\big\}\Big)\bigg)
E_{n-p+1}(f_\beta^\psi)_{L_1},
$$
де
$K_{q,p}(s)$ і $\sigma(s,p)$ визначаються формулами {\rm(\ref{S4.T.1.3})} і {\rm(\ref{S4.T.1.4})} відповідно, а $O(1)$ --- величина, рівномірно обмежена відносно всіх розглядуваних параметрів.

{\textbf{\textit{Наслідок {\rm\bf6}.}}}  \it Нехай множини $L^\psi_\beta L_1, \ \beta\in \mathbb{R}$ породжуються коефіцієнтами $\psi(k)=\frac{q^k}{k}$ ядер $N_{q,\beta}(t)$ вигляду {\rm(\ref{S4.T.1.6(11)})}. Тоді для довільних ${f\in L^\psi_\beta L_1}$, $n,p\in\mathbb{N}$, $p\leq n, \  1\leq s\leq\infty$ справедлива нерівність
$$
\|\rho_{n,p}(f;x)\|_{L_s}\leq$$$$\leq\frac{q^{n-p+1}}{p(n-p+1)}\big(\frac{\|\cos t\|_{s}}{\pi^{1+1/s}}K_{q,p}(s)+
O(1)\big(\frac{q}{(n-p+1)(1-q)^{\sigma(s,p)}}+$$$$+\frac{q}{(n-p+1)(1-q)^2}\min\big\{p,\frac{1}{1-q}\big\}\big)\big)
E_{n-p+1}(f_\beta^\psi)_{L_1},
$$
де
$K_{q,p}(s)$ і $\sigma(s,p)$  визначаються формулами {\rm(\ref{S4.T.1.3})} і {\rm(\ref{S4.T.1.4})} відповідно, а $O(1)$ --- величина, рівномірно обмежена відносно всіх розглядуваних параметрів.

\bf Теорема 3. \it{ Нехай $\psi\in {\cal D}_0$, $\beta\in\mathbb{R}$. Тоді для довільних $f\in C_\beta^\psi L_s$, ${1\leq s<\infty}$, $n,p\in\mathbb{N}$, $p\leq n$ справедлива нерівність
\begin{equation}\label{S4.T.2.1}
\|\rho_{n,p}(f;x)\|_C\leq$$$$\leq\Big(\frac{\|\cos t\|_{s'}}{\pi p}\psi (n-p+1)+O(1)\sum\limits_{k=n-p+2}^\infty\psi(k)\tau_{n,p}(k)\Big)E_{n-p+1}(f_\beta^\psi)_{L_s}.
\end{equation}

При цьому для будь-якої функції $f\in C^\psi_\beta L_s,  1\leq s<\infty$ і довільних $n,p\in \mathbb{N}$, $p\leq n$ в множині
$C^\psi_\beta L_s,  1\leq s<\infty$, знайдеться функція ${F(x)=F(f;n;p;x)}$ така, що
$E_{n-p+1}{(F_\beta^\psi)}_{L_s}=E_{n-p+1}{(f_\beta^\psi)}_{L_s}$, і для неї при $n-p\rightarrow\infty$ виконується рівність
\begin{equation}\label{S4.T.2.2}
\|\rho_{n,p}(F;x)\|_C=$$$$=\Big(\frac{\|\cos t\|_{s'}}{\pi p}\psi (n-p+1)+O(1)\sum\limits_{k=n-p+2}^\infty\psi(k)\tau_{n,p}(k)\Big)E_{n-p+1}(F_\beta^\psi)_{L_s},
\end{equation}
де $s'=\frac{s}{s-1}$.
У {\rm(\ref{S4.T.2.1})} і {\rm(\ref{S4.T.2.2})} коефіцієнти $\tau_{n,p}(k)$ визначаються рівністю {\rm(\ref{S4.T.1.9})}, а $O(1)$ --- величини, рівномірно обмежені відносно всіх розглядуваних параметрів.

{\textbf{\textit{Доведення теореми {\rm\bf3}.}}} \ \rm Перепишемо (\ref{S4.T.1.7}) у вигляді
\begin{equation}\label{S4.T.2.4}
\rho_{n,p}(f;x)=\frac{\psi(n-p+1)}{\pi p}\int\limits_{-\pi }^{\pi }f_\beta^\psi(x-t)\cos\Big((n-p+1)t-\frac{\beta\pi}{2}\Big)dt+$$$$+
\frac{1}{\pi}\int\limits_{-\pi }^{\pi }f_\beta^\psi(x-t)\Psi_{2,n,p}(t)dt,
\end{equation}
де $\Psi_{2,n,p}(t)$ означається рівністю (\ref{S4.T.1.8}).

Функції $\cos\big((n-p+1)t-\frac{\beta\pi}{2}\big)$ та $\Psi_{2,n,p}(t)$ ортогональні до будь-якого тригонометричного полінома $t_{n-p}$ порядку не вищого  $n-p$, тому в силу $(\ref{S4.T.2.4})$
\begin{equation}\label{S4.T.2.5}
\rho_{n,p}(f;x)=\frac{\psi(n-p+1)}{\pi p}\int\limits_{-\pi }^{\pi }\delta_{n,p}(x-t)\cos\Big((n-p+1)t-\frac{\beta\pi}{2}\Big)dt+$$$$+
\frac{1}{\pi}\int\limits_{-\pi }^{\pi }\delta_{n,p}(x-t)\Psi_{2,n,p}(t)dt,
\end{equation}
де $\delta_{n,p}(\cdot)$ визначається рівністю (\ref{S4.T.1.16}).

Обравши в (\ref{S4.T.2.5}) у ролі $t_{n-p}(\cdot)$ поліном $t_{n-p}^\ast(\cdot)$ найкращого наближення в просторі $L_s$ функції $f_\beta^\psi(\cdot)$ та використовуючи формули (\ref{S4.T.1.9}) і (\ref{S4.T.1.18'}),  інтеграли рівності (\ref{S4.T.2.5}) оцінимо наступним чином:
\begin{equation}\label{S4.T.2.9}
\big\|\int\limits_{-\pi }^{\pi }\delta_{n,p}(x-t)\cos(n-p+1)t\, dt\big\|_C\leq\|\delta_{n,p}\|_s\|\cos t\|_{s'}\leq\|\cos t\|_{s'}E_{n-p+1}(f_\beta^\psi)_{L_s},
\end{equation}
\begin{equation}\label{S4.T.2.7}
\big\|\frac{1}{\pi}\int\limits_{-\pi }^{\pi }\delta_{n,p}(x-t)\Psi_{2,n,p}(t)dt\big
\|_C\leq$$$$\leq \frac{1}{\pi}\|\delta_{n,p}\|_s\|\Psi_{2,n,p}\|_{s'}\leq \frac{2^{1/s'}}{\pi^{1/s}}\sum\limits_{k=n-p+2}^{\infty}\psi(k)\tau_{n,p}(k)E_{n-p+1}(f_\beta^\psi)_{L_s}.
\end{equation}
Об'єднуючи (\ref{S4.T.2.9}) і (\ref{S4.T.2.7}) отримуємо (\ref{S4.T.2.1}).

Доведемо другу частину теореми. З інтегрального зображення (\ref{S4.T.2.4}) і факту ортогональності функції $\Psi_{2,n,p}(t)$ до будь-якого тригонометричного полінома $t_{n-p}\in\mathcal{T}_{2(n-p)+1}$ випливає, що для довільної функції $f\in C_\beta^\psi L_s, \ \ 1\leq s<\infty$, $\psi\in{\cal D}_0$, $\beta\in\mathbb{R}$ виконується рівність
\begin{equation}\label{S4.T.2.10}
|\rho_{n,p}(f;x)|=\frac{\psi(n-p+1)}{\pi p}\Big|\int\limits_{-\pi }^{\pi }f_\beta^\psi(x-t)\cos\big((n-p+1)t-\frac{\beta\pi}{2}\big)dt\Big|+$$$$+
O(1)\sum\limits_{k=n-p+2}^\infty\psi(k)\tau_{n,p}(k)E_{n-p+1}(f_\beta^\psi)_{L_s}.
\end{equation}
Враховуючи (\ref{S4.T.2.10}), щоб переконатися в справедливості (\ref{S4.T.2.2}) досить показати, що якою б не була функція $\varphi\in L_s, \ \ 1\leq s<\infty$ знайдеться функція ${\Phi(\cdot)=\Phi(\varphi;\cdot)}$, для якої при всіх  $n,p\in\mathbb{N}$, $p\leq n$
\begin{equation}\label{S4.T.2.11}
E_{n-p+1}(\Phi)_{L_s}=E_{n-p+1}(\varphi)_{L_s}
\end{equation}
і, крім того, має місце рівність
\begin{equation}\label{S4.T.2.12}
\big|\int\limits_{-\pi }^{\pi }\Phi(t)\cos\Big((n-p+1)t+\frac{\beta\pi}{2}\big)dt\Big|=\|\cos t\|_{s'}E_{n-p+1}(\varphi)_{L_s},
\end{equation}

В якості $\Phi(\cdot)$  розглянемо функцію
\begin{equation}\label{S4.T.2.13}
\Phi(t)=$$$$=\|\cos t\|_{s'}^{1-s'}\big|\cos\big((n-p+1)t+\frac{\beta\pi}{2}\big)\big|^{s'-1}
\mathrm{sign}\cos\big((n-p+1)t+\frac{\beta\pi}{2}\big)E_{n-p+1}(\varphi)_{L_s}.
\end{equation}
Для неї
\begin{equation}\label{S4.T.2.14}
\|\Phi(t)\|_s=\|\cos t\|_{s'}^{1-s'}\Big(\int\limits_{-\pi }^{\pi }\big|\cos\big((n-p+1)t+\frac{\beta\pi}{2}\Big)\big|^{(s'-1)s}dt\big)^\frac{1}{s}E_{n-p+1}(\varphi)_{L_s}=$$$$
=\|\cos t\|_{s'}^{1-s'}\big\|\cos\big((n-p+1)t+\frac{\beta\pi}{2}\big)\big\|^{s'-1}_{s'}E_{n-p+1}(\varphi)_{L_s}=E_{n-p+1}(\varphi)_{L_s}.
\end{equation}
Крім того, оскільки для довільного $t_{n-p}\in\mathcal{T}_{2(n-p)+1}$
$$
\int\limits_{-\pi}^{\pi}t_{n-p}(\tau)|\Phi(\tau)|^{s-1}\mathrm{sign}\Phi(\tau)d\tau=$$$$=
\Big(\|\cos t\|_{s'}^{1-s'}E_{n-p+1}(\varphi)_{L_s}\Big)^{s-1}\int\limits_{-\pi}^{\pi}t_{n-p}(\tau)\cos\Big((n-p+1)\tau+\frac{\beta\pi}{2}\Big)d\tau=0,
$$
то на підставі теореми 1.4.5 роботи [\ref{Korn}, c. 28] можемо зробити висновок, що поліном ${t_{n-p}^\ast\equiv0}$ є
поліномом найкращого наближення функції $\Phi(t)$ в метриці простору ${L_s, 1\leq s<\infty}$.
Отже,  з урахуванням (\ref{S4.T.2.14}),
\begin{equation}\label{S4.T.2.14'}
E_{n-p+1}(\Phi)_{L_s}=\|\Phi\|_s=E_{n-p+1}(\varphi)_{L_s}.
\end{equation}
В силу (\ref{S4.T.2.13}),
$$
\Big|\int\limits_{-\pi}^{\pi}\Phi(t)\cos\big((n-p+1)t+\frac{\beta\pi}{2}\big)dt\Big|=$$$$=
\|\cos t\|_{s'}^{1-s'}E_{n-p+1}(\varphi)_{L_s}\Big|\int\limits_{-\pi}^{\pi}\big|\cos\big((n-p+1)t+\frac{\beta\pi}{2}\big)\big|^{s'-1}
\times$$$$\times\mathrm{sign}\cos\big((n-p+1)t+\frac{\beta\pi}{2}\big)\cos\big((n-p+1)t+\frac{\beta\pi}{2}\big)dt\Big|=$$$$=
\|\cos t\|_{s'}^{1-s'}E_{n-p+1}(\varphi)_{L_s}\int\limits_{-\pi}^{\pi}\big|\cos\big((n-p+1)t+\frac{\beta\pi}{2}\big)\big|^{s'}dt=\|\cos t\|_{s'}E_{n-p+1}(\varphi)_{L_s},
$$
З останніх  співвідношень випливає (\ref{S4.T.2.12}). Теорему 3 доведено.

\bf Теорема 4. \it{ Нехай $\psi\in {\cal D}_0$, $\beta\in\mathbb{R}$. Тоді для довільних $f\in C_\beta^\psi L_\infty$, $n,p\in\mathbb{N}$, $p\leq n$  справедлива нерівність
\begin{equation}\label{S4.T.3.1(1)}
\|\rho_{n,p}(f;x)\|_C\leq$$$$\leq\frac{1}{p}\bigg(\frac{4}{\pi}\psi(n-p+1)+O(1)\Big(\frac{\psi^2(n-p+2)}{\psi(n-p+1)}
+p\!\!\!\!\sum\limits_{k=n-p+3}^\infty\psi(k)\tau_{n,p}(k)\Big)\bigg)E_{n-p+1}(f_\beta^\psi)_{L_\infty}.
\end{equation}

При цьому для будь-якої функції $f\in C_\beta^\psi L_\infty$ і довільних ${n,p\in \mathbb{N}}$, $p\leq n$  знайдеться функція ${F(x)=F(f;n;p;x)}\in C_\beta^\psi C$ така, що ${E_{n-p+1}{(F_\beta^\psi)}_C=E_{n-p+1}{(f_\beta^\psi)}_{L_\infty}}$ і для неї при $n-p\rightarrow\infty$ виконується рівність
\begin{equation}\label{S4.T.3.2(1)}
\|\rho_{n,p}(F;x)\|_C=$$$$=\frac{1}{p}\bigg(\frac{4}{\pi}\psi(n-p+1)+O(1)\Big(\frac{\psi^2(n-p+2)}{\psi(n-p+1)}
+p\!\!\!\!\sum\limits_{k=n-p+3}^\infty\psi(k)\tau_{n,p}(k)\Big)\bigg)E_{n-p+1}(F_\beta^\psi)_C.
\end{equation}

У {\rm(\ref{S4.T.3.1(1)})} і {\rm(\ref{S4.T.3.2(1)})} коефіцієнти $\tau_{n,p}(k)$ визначаються рівністю {\rm(\ref{S4.T.1.9})},
а $O(1)$ --- величини, рівномірно обмежені відносно всіх розглядуваних параметрів.

{\textbf{\textit{Доведення теореми {\rm\bf4}.}}} \ \rm Нехай $f\in C_\beta^\psi L_\infty, \ \psi\in{\cal D}_0$. Виходячи з (\ref{S4.T.1.7}), та враховуючи факт ортогональності функції $\Psi_{1,n,p}(t)$ до будь-якого полінома $t_{n-p}$ порядку не вищого за $n-p$, можемо записати
\begin{equation}\label{S4.T.3.3}
\rho_{n,p}(f;x)=\frac{1}{\pi}\int\limits_{-\pi}^\pi \delta_{n,p}(x-t)\Psi_{1,n,p}(t)dt=$$$$=\frac{1}{\pi p}\int\limits_{-\pi}^\pi\delta_{n,p}(x-t)\big(\psi(n-p+1)
\cos\big((n-p+1)t-\frac{\beta\pi}{2}\big)+$$$$+2\psi(n-p+2)\cos\big((n-p+2)t-\frac{\beta\pi}{2}\big)+p\!\!\!\!\sum\limits_{k=n-p+3}^\infty\!\!\!\tau_{n,p}(k)\psi(k)\cos\big(kt-\frac{\beta\pi}{2}\big)\big)dt,
\end{equation}
де $\tau_{n,p}(k)$ і $\delta_{n,p}(\cdot)$ визначаються рівностями (\ref{S4.T.1.9}) і (\ref{S4.T.1.16}) відповідно.

Обравши в (\ref{S4.T.3.3}) у ролі $t_{n-p}$ поліном $t_{n-p}^\ast$ найкращого наближення в просторі $L_\infty$ функції $f_\beta^\psi$ і застосувавши нерівність (\ref{S4.T.1.18'}) при $s=\infty$, отримуємо оцінку
\begin{equation}\label{S4.T.3.7}
\|\rho_{n,p}(f;x)\|_C\leq$$$$\leq\frac{1}{\pi p}\big(\big\|\psi(n-p+1)
\cos\big((n-p+1)t-\frac{\beta\pi}{2}\big)+2\psi(n-p+2)\cos\big((n-p+2)t-\frac{\beta\pi}{2}\big)\big\|_1+$$$$+p\big\|\sum\limits_{k=n-p+3}^\infty\tau_{n,p}(k)\psi(k)\cos\big(kt-\frac{\beta\pi}{2}\big)\big\|_1\big)E_{n-p+1}{(f_\beta^\psi)}_{L_\infty}\leq$$$$
\leq\frac{1}{\pi p}\Big(\big\|\psi(n-p+1)
\cos(n-p+1)t+$$$$+2\psi(n-p+2)\cos\big((n-p+2)t+\frac{\beta\pi}{2(n-p+1)}\big)\big\|_1+$$$$+O(1)p\!\!\!\!\sum\limits_{k=n-p+3}^\infty\tau_{n,p}(k)\psi(k)\Big)E_{n-p+1}{(f_\beta^\psi)}_{L_\infty},
\end{equation}
Як випливає з роботи С.О.~Теляковського [\ref{Teliacovs 1989}, c. 512-513],
\begin{equation}\label{S4.T.3.8}
\big\|\psi(n-p+1)\cos(n-p+1)t+2\psi(n-p+2)\cos\big((n-p+2)t+\frac{\beta\pi}{2(n-p+1)}\big)\big\|_1+$$$$+O(1)p\!\!\!\!\sum\limits_{k=n-p+3}^\infty\tau_{n,p}(k)\psi(k)\leq$$$$\leq 4\psi(n-p+1)+O(1)\Big(\frac{\psi^2(n-p+2)}{\psi(n-p+1)}+p\!\!\!\!\sum\limits_{k=n-p+3}^\infty\!\!\!\tau_{n,p}(k)\psi(k)\Big).
\end{equation}
Співвідношення (\ref{S4.T.3.7}) і (\ref{S4.T.3.8}) доводять нерівність (\ref{S4.T.3.1(1)}).

Доведемо другу частину теореми 4. Виходячи з інтегрального зображення (\ref{S4.T.3.3}) і використовуючи факт ортогональності функції ${\sum\limits_{k=n-p+3}^\infty\!\!\!\tau_{n,p}(k)\psi(k)\cos\big(kt-\frac{\beta\pi}{2}\big)}$ до будь-якого тригонометричного полінома $t_{n-p}$ порядку не вищого $n-p$, для довільної функції $f$  з множини $C_\beta^\psi L_\infty$, $\psi\in {\cal D}_0$, $\beta\in\mathbb{R}$ виконується рівність
\begin{equation}\label{S4.T.3.9}
|\rho_{n,p}(f;x)|=\frac{\psi(n-p+1)}{\pi p}\Big|\int\limits_{-\pi}^\pi f_\beta^\psi(x-t)\big(
\cos\big((n-p+1)t-\frac{\beta\pi}{2}\big)+$$$$+2\frac{\psi(n-p+2)}{\psi(n-p+1)}\cos\big((n-p+2)t-\frac{\beta\pi}{2}\big)\big)dt\Big|+
O(1)\!\!\!\sum\limits_{k=n-p+3}^\infty\!\!\!\tau_{n,p}(k)\psi(k)E_{n-p+1}{(f_\beta^\psi)}_{L_\infty}.
\end{equation}
Для доведення (\ref{S4.T.3.2(1)}), з урахуванням (\ref{S4.T.3.9}), досить встановити, що для довільної ${\varphi\in L_\infty^0=\{\varphi\in L_\infty: \varphi\bot1\}}$ існує функція $\Phi(\cdot)=\Phi(\varphi;\cdot)\in C$ для якої при всіх $n,p\in \mathbb{N}$, $p\leq n$
$$
E_{n-p+1}(\Phi)_C=E_{n-p+1}(\varphi)_{L_\infty}
$$
і, крім того, при $n-p\rightarrow\infty$ має місце рівність
\begin{equation}\label{S4.T.3.10}
\big|\int\limits_{-\pi}^\pi \Phi(t)\big(
\cos\big((n-p+1)t+\frac{\beta\pi}{2}\big)+2\frac{\psi(n-p+2)}{\psi(n-p+1)}\cos\big((n-p+2)t+\frac{\beta\pi}{2}\big)\big)dt\big|=$$$$=
\bigg(4+O(1)\Big(\frac{\psi(n-p+2)}{\psi(n-p+1)}\Big)^2\bigg)E_{n-p+1}(\varphi)_{L_\infty}.
\end{equation}

Покладемо
$$
\varphi_0(t)=\mathrm{sign}\cos\Big((n-p+1)t+\frac{\beta\pi}{2}\Big)E_{n-p+1}(\varphi)_\infty
$$
і через $\varphi_\delta(t)$ позначимо $2\pi$-періодичну функцію, яка збігається з $\varphi_0(t)$ скрізь, за виключенням $\delta$-околів (${0<\delta<\frac{\pi}{2(n-p+1)}}$)  точок $t_k=\frac{(2k+1-\beta)\pi}{2(n-p+1)}, \ \ k\in\mathbb{Z}$, де вона лінійна і її графік сполучає точки $(t_k-\delta, \varphi_0(t_k-\delta))$ і $(t_k+\delta, \varphi_0(t_k+\delta))$. Функція $\varphi_\delta(t)$ неперервна і у точках ${\tau_k=\frac{(2k-\beta)\pi}{2(n-p+1)}, \ \ k=1,2,...,2(n-p+1)}$ періоду $\big(-\frac{\beta\pi}{2(n-p+1)}, 2\pi-\frac{\beta\pi}{2(n-p+1)}\big]$ досягає по абсолютній величині максимального значення, яке дорівнює $E_{n-p+1}(\varphi)_\infty$, почергово змінюючи знак. Тому її поліном найкращого рівномірного наближення порядку не вищого $n-p$, згідно з критерієм Чебишова, є поліном, що тотожно дорівнює нулю і, отже,
\begin{equation}\label{S4.T.3.11}
E_{n-p+1}(\varphi_\delta)_C=\|\varphi_\delta\|_C=E_{n-p+1}(\varphi)_\infty.
\end{equation}
Враховуючи (\ref{S4.T.1.18'}) і (\ref{S4.T.3.11}), одержуємо
\begin{equation}\label{S4.T.3.12}
\big|\int\limits_{-\pi}^\pi \varphi_\delta(t)\big(
\cos\big((n-p+1)t+\frac{\beta\pi}{2}\big)+2\frac{\psi(n-p+2)}{\psi(n-p+1)}\cos\big((n-p+2)t+\frac{\beta\pi}{2}\big)\big)dt\big|\leq$$$$\leq
\int\limits_{-\pi}^\pi \big|\cos(n-p+1)t+2\frac{\psi(n-p+2)}{\psi(n-p+1)}\cos\big((n-p+2)t-\frac{\beta\pi}{2(n-p+1)}\big)\big|dt \times$$$$\times E_{n-p+1}(\varphi)_{L_\infty}.
\end{equation}

Із нерівності (19) роботи [\ref{Teliacovs 1989}, c. 513], випливає оцінка
\begin{equation}\label{S4.T.3.13}
\int\limits_{-\pi}^\pi \big|\cos(n-p+1)t+2\frac{\psi(n-p+2)}{\psi(n-p+1)}\cos\big((n-p+2)t-\frac{\beta\pi}{2(n-p+1)}\big)\big|dt\leq$$$$\leq 4+O(1)\Big(\frac{\psi(n-p+2)}{\psi(n-p+1)}\Big)^2.
\end{equation}
З іншого боку,
\begin{equation}\label{S4.T.3.14}
\big|\int\limits_{-\pi}^\pi \varphi_\delta(t)\big(
\cos\big((n-p+1)t+\frac{\beta\pi}{2}\big)+2\frac{\psi(n-p+2)}{\psi(n-p+1)}\cos\big((n-p+2)t+\frac{\beta\pi}{2}\big)\big)dt\big|=$$$$
=\big|\int\limits_{-\pi}^\pi \varphi_0(t)\big(
\cos\big((n-p+1)t+\frac{\beta\pi}{2}\big)+2\frac{\psi(n-p+2)}{\psi(n-p+1)}\cos\big((n-p+2)t+\frac{\beta\pi}{2}\big)\big)dt\big|+$$$$+O(1)r_{n,p}(\delta),
\end{equation}
де
\begin{equation}\label{S4.T.3.15}
r_{n,p}(\delta)=\big|\int\limits_{-\pi}^\pi (\varphi_\delta(t)-\varphi_0(t))\big(
\cos\big((n-p+1)t+\frac{\beta\pi}{2}\big)+$$$$+2\frac{\psi(n-p+2)}{\psi(n-p+1)}\cos\big((n-p+2)t+\frac{\beta\pi}{2}\big)\big)dt\big|.
\end{equation}
Оскільки $\psi\in {\cal D}_0$, то для досить великих номерів $n-p$ справджується нерівність  $\frac{\psi(n-p+2)}{\psi(n-p+1)}<1$ і, отже,
\begin{equation}\label{S4.T.3.16}
r_{n,p}(\delta)<3\int\limits_{-\pi}^{\pi}|\varphi_\delta(t)-\varphi_0(t)|dt\leq6(n-p+1)\delta E_{n-p+1}(\varphi)_\infty.
\end{equation}
Вибравши $\delta$ настільки малим, щоб виконувалась умова
\begin{equation}\label{S4.T.3.17}
0<\delta<\frac{1}{n-p+1}\Big(\frac{\psi(n-p+2)}{\psi(n-p+1)}\Big)^2,
\end{equation}
із (\ref{S4.T.3.16}) одержимо оцінку
\begin{equation}\label{S4.T.3.18}
r_{n,p}(\delta)=O(1)\Big(\frac{\psi(n-p+2)}{\psi(n-p+1)}\Big)^2E_{n-p+1}(\varphi)_\infty.
\end{equation}
Оскільки $\int\limits_{-\pi}^{\pi}\varphi_0(t)\cos\big((n-p+2)t+\frac{\beta\pi}{2}\big)dt=0$, то
\begin{equation}\label{S4.T.3.19}
\big|\int\limits_{-\pi}^\pi \varphi_0(t)\big(
\cos\big((n-p+1)t+\frac{\beta\pi}{2}\big)+2\frac{\psi(n-p+2)}{\psi(n-p+1)}\cos\big((n-p+2)t+\frac{\beta\pi}{2}\big)\big)dt\big|=$$$$
=\big|\int\limits_{-\pi}^{\pi}\varphi_0(t)\cos\big((n-p+1)t+\frac{\beta\pi}{2}\big)dt\big|=$$$$=
\int\limits_{-\pi}^{\pi}\big|\cos\big((n-p+1)t+\frac{\beta\pi}{2}\big)\big|dtE_{n-p+1}(\varphi)_\infty=4E_{n-p+1}(\varphi)_\infty.
\end{equation}
Із формул (\ref{S4.T.3.12})--(\ref{S4.T.3.14}), (\ref{S4.T.3.18}) і (\ref{S4.T.3.19}) бачимо, що для функції ${\Phi(t)=\varphi_\delta(t)}$ у якій параметр $\delta$ задовольняє умову (\ref{S4.T.3.17}), при $n-p+1\rightarrow\infty$ має місце рівність (\ref{S4.T.3.10}), а, отже, і (\ref{S4.T.3.2(1)}). Теорему 4 доведено.

\rm Оскільки для сум $\sum\limits_{k=n-p+j}^\infty\tau_{n,p}(k)\psi(k)$, які фігурують в теоремах {\rm3 та 4,}  мають місце рівності
\begin{equation}\label{S4.T.3.24'}
   \sum\limits_{k=n-p+j}^{\infty}\tau_{n,p}(k)\psi(k)=\left\{\begin{array}{ll}
   \sum\limits_{k=n-p+j}^{n-1}\frac{k-n+p}{p}\psi (k)+\sum\limits_{k=n}^{\infty}\psi (k), & \ p>j, \ \ j\in\mathbb{N},\\
   \sum\limits_{k=n-p+j}^{\infty}\psi(k), & \ p\leq j, \ \ j\in\mathbb{N},
  \end{array}\right.
\end{equation}
то, як неважко переконатися, для них справедлива наступна оцінка зверху:
$$
\sum\limits_{k=n-p+j}^{\infty}\tau_{n,p}(k)\psi(k)\leqslant \min\big\{\sum\limits_{k=n-p+j}^{\infty}\!\!\!\psi (k),\frac{1}{p}\sum\limits_{k=n-p+j}^{\infty}(k-n+p)\psi (k)\big\},\ \ \ \ j\in \mathbb{N}.
$$
Отже, в співвідношеннях (\ref{S4.T.2.1}) і (\ref{S4.T.2.2}) теореми 3 та співвідношеннях (\ref{S4.T.3.1(1)}) і (\ref{S4.T.3.2(1)}) теореми 4 величини $O(1)\!\!\!\sum\limits_{k=n-p+j}^\infty\tau_{n,p}(k)\psi(k)$ можна замінити на ${O(1)\min\big\{\sum\limits_{k=n-p+j}^{\infty}\!\!\!\psi (k),\frac{1}{p}\sum\limits_{k=n-p+j}^{\infty}(k-n+p)\psi (k)\big\}}$,  $j=2,3$.

Нерівності (\ref{S4.T.2.1}) і (\ref{S4.T.3.1(1)}) залишаються асимптотично непокращуваними не тільки на усіх множинах  $C_\beta^\psi L_s, \ 1\leq s<\infty$ та $C_\beta^\psi L_\infty$  при $\psi\in {\cal D}_0$, але і на таких важливих їх підмножинах, як класи $C_{\beta,s}^\psi$ і $C_{\beta,\infty}^\psi$. Це випливає із наступних міркувань. Розглянемо точні верхні межі в обох частинах нерівності (\ref{S4.T.2.1}) по класу ${C_{\beta,s}^\psi \ 1\leq s<\infty}$ і точні верхні межі в обох частинах нерівності (\ref{S4.T.3.1(1)}) по класу $C_{\beta,\infty}^\psi$. В результаті одержуємо нерівності
\begin{equation}\label{S4.T.3.20}
{\cal E}(C^\psi_{\beta,s};V_{n,p})_C\leq\bigg(\frac{\|\cos t\|_{s'}}{\pi p}\psi (n-p+1)+O(1)\!\!\!\sum\limits_{k=n-p+2}^\infty\!\!\!\psi(k)\tau_{n,p}(k)\bigg), \ 1\leq s<\infty,
\end{equation}
\begin{equation}\label{S4.T.3.21}
{\cal E}(C^\psi_{\beta,\infty};V_{n,p})_C\leq\frac{1}{p}\bigg(\frac{4}{\pi}\psi(n-p+1)+O(1)\Big(\frac{\psi^2(n-p+2)}{\psi(n-p+1)}
+p\!\!\!\!\sum\limits_{k=n-p+3}^\infty\!\!\!\psi(k)\tau_{n,p}(k)\Big)\bigg).
\end{equation}
Зіставляючи два останні співвідношення відповідно з рівностями (24) і (7) роботи [\ref{Serduk Ovsii 2008},~c.~337, 341], приходимо до висновку, що у співвідношеннях  (\ref{S4.T.3.20}) і (\ref{S4.T.3.21}) можна поставити знак "дорівнює".

У випадку, коли $\psi (k)=e^{-\alpha k^r},$ $\alpha >0$, $r>1$ (тобто, коли $C_\beta^\psi L_s=C_\beta^{\alpha,r} L_s$), як показано в [\ref{Serduk Ovsii 2008}, c. 348], для довільних $ \alpha>0,\ r>1$ має місце оцінка
\begin{equation}\label{S4.T.3.22}
\sum\limits_{k=n-p+j}^{\infty}\!\!\!\!\psi(k)\tau_{n,p}(k)=\sum\limits_{k=n-p+j}^{\infty}\!\!\!\!e^{-\alpha k^r}\tau_{n,p}(k)=$$$$=O(1)j^{\sigma(p)-1}\big(1+\frac{1}{\alpha r(n-p+j)^{r-1}}\big)^{\sigma (p)}e^{-\alpha (n-p+j)^r}, \ j\in\mathbb{N},
\end{equation}
де
\begin{equation}\label{S4.T.3.24}
  \sigma (p)\stackrel{\mathrm{df}}{=}\left\{\begin{array}{ll}
    1, & \mbox{якщо},\ p\leqslant j, \\
    2, & \mbox{якщо},\ p>j,\ \ \ j\in \mathbb{N},
  \end{array}\right.
\end{equation}
а коефіцієнти $\tau_{n,p}(k)$ означаються формулою (\ref{S4.T.1.9}). Отже, в силу (\ref{S4.T.3.22}) із теорем 3 та 4 одержуємо наступні твердження.

{\textbf{\textit{Наслідок {\rm\bf7}.}}} \it {Нехай $\alpha >0,\ r>1,$ $\beta\in \mathbb{R},\ p,n\in \mathbb{N},\ p\leq n$ і $1\leq s<\infty.$ Тоді
для довільних $f\in C_{\beta}^{\alpha,r}L_s$ справедлива нерівність
\begin{equation}\label{S4.T.3.25}
\|\rho_{n,p}(f;x)\|_C\leq \frac{e^{-\alpha (n-p+1)^r}}{p}\bigg(\frac{\|\cos t\|_{s'}}{\pi}+$$$$+O(1)\Big(1+\frac{1}{\alpha r(n-p+2)^{r-1}}\Big)^{\sigma (p)}\frac{e^{\alpha (n-p+1)^r}}{e^{\alpha (n-p+2)^r}}\bigg)E_{n-p+1}(f_\beta^{\alpha,r})_{L_s},
\end{equation}
де $\sigma (p)$ означається формулою {\rm(\ref{S4.T.3.24})}, $s'=s/(s-1)$.

При цьому для будь-якої функції $f\in C_{\beta}^{\alpha,r}L_s, \ \alpha >0,\ r>1$,  $1\leq s<\infty$ і довільних $n,p\in \mathbb{N}$, $p\leq n$ в множині $C_{\beta}^{\alpha,r}L_s,  1\leq s<\infty$, знайдеться функція ${F(x)=F(f;n;p;x)}$ така, що $E_{n-p+1}{(F_\beta^{\alpha,r})}_{L_s}=E_{n-p+1}{(f_\beta^{\alpha,r})}_{L_s}$, і для неї при ${n-p\rightarrow\infty}$ виконується рівність
\begin{equation}\label{S4.T.3.25'}
\|\rho_{n,p}(F;x)\|_C=\frac{e^{-\alpha (n-p+1)^r}}{p}\bigg(\frac{\|\cos t\|_{s'}}{\pi}+$$$$+O(1)\Big(1+\frac{1}{\alpha r(n-p+2)^{r-1}}\Big)^{\sigma (p)}\frac{e^{\alpha (n-p+1)^r}}{e^{\alpha (n-p+2)^r}}\bigg)E_{n-p+1}(F_\beta^{\alpha,r})_{L_s}.
\end{equation}
У {\rm(\ref{S4.T.3.25})} і {\rm(\ref{S4.T.3.25'})} $O(1)$ --- величини, рівномірно обмежені відносно всіх розглядуваних параметрів.}

{\textbf{\textit{Наслідок {\rm\bf8}.}}}  \it Нехай $\alpha >0,\ r>1,$ $\beta\in \mathbb{R},\ p,n\in \mathbb{N},\ p\leq n$. Тоді
для довільних $f\in C_{\beta}^{\alpha,r}L_\infty$ справедлива нерівність
\begin{equation}\label{S4.T.3.26}
\|\rho_{n,p}(f;x)\|_C\leq \frac{e^{-\alpha (n-p+1)^r}}{p }\bigg(\frac{4}{\pi }+O(1)\Big(\frac{e^{2\alpha (n-p+1)^r}}{e^{2\alpha (n-p+2)^r}}+$$$$+
\big(1+\frac{1}{\alpha r(n-p+3)^{r-1}}\big)^{\sigma (p)}\frac{e^{\alpha (n-p+1)^r}}{e^{\alpha(n-p+3)^r}}\Big)\bigg)E_{n-p+1}(f_\beta^{\alpha,r})_{L_\infty},
\end{equation}
де $\sigma (p)$ означається формулою {\rm(\ref{S4.T.3.24})}.

При цьому для будь-якої функції $f\in C_{\beta}^{\alpha,r}L_\infty, \ \alpha >0,\ r>1$  і довільних ${n,p\in \mathbb{N}}$, $p\leq n$ в множині
$C_{\beta}^{\alpha,r}L_\infty$, знайдеться функція ${F(x)=F(f;n;p;x)}$ така, що $E_{n-p+1}{(F_\beta^{\alpha,r})}_{L_\infty}=E_{n-p+1}{(f_\beta^{\alpha,r})}_{L_\infty}$, і для неї при ${n-p\rightarrow\infty}$ виконується рівність
\begin{equation}\label{S4.T.3.27}
\|\rho_{n,p}(F;x)\|_C=\frac{e^{-\alpha (n-p+1)^r}}{p }\bigg(\frac{4}{\pi }+O(1)\Big(\frac{e^{2\alpha (n-p+1)^r}}{e^{2\alpha (n-p+2)^r}}+$$$$+
\big(1+\frac{1}{\alpha r(n-p+3)^{r-1}}\big)^{\sigma (p)}\frac{e^{\alpha (n-p+1)^r}}{e^{\alpha(n-p+3)^r}}\Big)\bigg)E_{n-p+1}(F_\beta^{\alpha,r})_{L_\infty}.
\end{equation}
У {\rm(\ref{S4.T.3.26})} і {\rm(\ref{S4.T.3.27})} $O(1)$ --- величини, рівномірно обмежені відносно всіх розглядуваних параметрів.} \rm

\newpage

\begin{enumerate}

\centerline {\bf Література}
\item \label{Step monog 1987} {\it  Степанец А.И.} Классификация и приближение периодических функций.~--- Киев: Наук. думка, 1987.~--- 268 c.

\item \label{Step monog 2002(1)}{\it Степанец А.И.\/} Методы теории приближений  В 2 ч. // Праці Ін-ту
математики НАН України. Т 40. --- К.: Ін-т математики НАН України, 2002.~--- Ч.I. --- 427 с.

\item \label{Falaleev 2001} {\it  Фалалеев Л.П.} О приближении  функций обобщенными операторами Абеля-Пуассона // Сиб. мат. журн. ---
2001. --- {\bf 42}, №4.~--- C.~926--936.

\item \label{Stech 1978} {\it Steckin S.B.} On the appoximation of
periodic functions by de la Vall\'{e}e Poussin sums // Anal. math. --- 1978. --- {\bf 4}. --- P. 61--74.

\item \label{Lebeg1910} {\it   Lebesgue H.} Sur la repr\'{e}sentation trigonom\'{e}trique approch\'{e}e des fonctions
satisfaisantes \'{a} une condition de Lipschitz // Bull. Soc. Math. France. -- 1910.~--- {\bf38}. --- P. 184--210.

\item \label{Valle Pyssen 1919} {\it Ch. de la Vall\'{e}e Poussin.} Lecons sur l'approximation des fonctions d'une variable r\'{e}elle. --- Paris: Gautier-Villars, 1919. --- 150 p.

\item \label{Nikol 1940} {\it Никольский С.М.} О некоторых методах приближения тригонометрическими суммами // Изв. АН СССР. Сер. мат.
--- 1940. --- {\bf 4}. --- С. 509--520.

\item \label{Stech 1951} {\it Стечкин С.Б.} О суммах Валле Пуссена // Докл. АН СССР. --- 1951. --- {\bf 80}.~---~С.~545--548.

\item \label{Gabis 1965} {\it Габисония О.Д.} О приближении функций многих переменных целыми функциями // Изв. вуз. Матем. --- 1965. ---
{\bf 2(45)}. --- С. 30--35.

\item \label{Zaharov 1968} {\it Захаров А.А.} Об оценке уклонения непрерывных периодических функций от сумм Валле Пуссена // Мат. заметки. ---
1968. --- {\bf 3}. --- C. 77--84.

\item \label{Oscolcov 1975} {\it Осколков К.И.} К неравенству Лебега в равномерной метрике и на множестве полной меры // Мат. заметки. ---
1975. --- {\bf 18}. --- C. 515--526.

\item \label{Stech 1961} {\it Стечкин С.Б.} О приближении периодических функций суммами Фейера // Тр. Матем. ин-та  АН
СССР. --- 1961. --- {\bf 62}. --- \ C.~48--60.

\item \label{Kolmogorov 1935 } {\it  Kolmogoroff A.N.} Zur Gr\"{o}ssenordnung des Restgliedes Fouriershen Reihen
differenzierbarer Funktionen // Ann. Math. --- 1935. --- {\bf 36}. --- S. 521--526.

\item \label{Nikol 1941} {\it Никольский С.М.} Асимптотическая оценка остатка при приближении суммами Фурье // Докл. АН СССР. ---
1941. --- {\bf 22}, №6.~--- C.~386--389.

\item \label{Nikol 1945} {\it Никольский С.М.} Приближение периодических функций тригонометрическими многочленами //
Тр. Матем. ин-та АН СССР. --- 1945. --- {\bf 15}. --- C.~1--~76.

\item \label{Timan 1951 } {\it Тиман А.Ф.} Обобщение некоторых результатов А.Н. Колмогорова и С.М.~Никольского //
Докл. АН СССР. --- 1951. --- {\bf 81}, №4. --- C.~509--511.

\item \label{Step,Rycasov monog 2007} {\it  Степанец А.И., Рукасов В.И., Чайченко С.О.} Приближения суммами Валле Пуссена // Праці Ін-ту
математики НАН України. --- 2007.--- {\bf68}. -- 386~с.

\item \label{RykasovHcay2002(1)} {\it Рукасов В.І., Чайченко С.О.} Наближення аналітичних періодичних функцій сумами Валле Пуссена//  Укр. мат. журн. -- 2002. -- {\bf 54}, №12.~-- C.~1653--1668.

\item \label{Step monog 2002(2)} {\it  Степанец А.И.} Методы теории приближений: В 2 ч. // Праці Ін-ту математики НАН України. -- 2002. -- {\bf40}. {\it ч}.II. -- 424 с.

\item \label{Rykasov 2003} {\it Рукасов В.И.} Приближение суммами Валле Пуссена классов аналитических функций / / Укр. мат. журн. -- 2003. -- {\bf55},  №6. -- С.~806--816.

\item \label{Serd 2004} {\it Сердюк А.С.}  Наближення інтегралів Пуассона сумами Валле Пуссена // Укр. мат. журн. -- 2004. -- {\bf 56}, №1.~-- C.~97--107.

\item \label{Serduk Ovsii 2008} {\it Сердюк А.С., Овсій Є.Ю.} Наближення на класах цілих функцій сумами Валле Пуссена //  Теорія наближення функцій та суміжні питання: Зб. праць Ін-ту математики НАН України. -- 2008. -- 5, №1. -- С.~334--351.

\item \label{Serd2009} {\it Сердюк А.С.}  Наближення інтегралів Пуассона сумами Валле Пуссена в рівномірній та інтегральних метриках
// Доп. НАН України -- 2009. -- №6.~-- C.~34--39.

\item \label{Serd 2010} {\it Сердюк А.С.} Приближение интегралов Пуассона суммами Валле Пуссена в равномерной и интегральных метриках  // Укр. мат. журн. -- 2010. -- {\bf 62}, №12.~-- C.~1672--1686.

\item \label{Serduk Ovsii 2011(1)} {\it  Serdyuk A.S., Ovsii Ie.Yu.} Uniform approximation of Poisson integrals of functions from the class $H_\omega$ by de la Vall\'{e}e Poussin sums // Arxiv preprint, arXiv:1104.3060, 2011. -- 19 p.

\item\label{Serd_Ovsii_2011}{\it Serdyuk A.S., Ovsii Ie.Yu., Musienko A.P.\/} Approximation of classes of analytic functions by de la Vallee
Poussin sums in uniform metric // Arxiv preprint, arXiv:1112.0967, 2011. -- 14 p.

\item \label{Stepan i Serduk 2000} {\it Степанец А.И., Сердюк А.С.}  Неравенства Лебега для интегралов Пуассона // Укр. мат. журн. -- 2000.
-- {\bf52}, №6. -- C. 798--808.

\item \label{Serduk Mysienko 2010} {\it Сердюк А.С., Мусієнко А.П.} Нерівності типу Лебега для сум Валле Пуссена при наближенні інтегралів Пуассона //  Теорія наближення функцій та суміжні питання: Зб. праць Ін-ту математики НАН України. -- 2010.~-- 7, №1. -- С.~298--316.

\item \label{Serd_Step 2000}{\it Степанец А.И., Сердюк А.С.}  Приближение суммами Фурье и наилучшие приближения на классах аналитических функций
// Укр. мат. журн.~-- 2000.~-- {\bf 52}, №3.~-- C.~375--395.

\item \label{Timan 2009}{\it  Тиман М.Ф.\/} Аппроксимация и свойства периодических функций. -- Киев: Наук. думка, 2009. -- 376 с.

\item \label{Axiezer monog}{\it  Ахиезер Н.И.\/} Лекции по теории аппроксимации. -- М.: Наука, 1965.~--  408~c.

\item \label{Serduk Cay 2011}{\it  Сердюк А.С., Чайченко С.О.\/} Наближення класів аналітичних функцій лінійним методом спеціального вигляду// Укр. мат. журн.-- 2011.-- {\bf 63}, №~1. -- C.~102--109.

\item \label{Serduk Sokolenko 2010}{\it  Serdyuk А.S., Sokolenko I.V.\/} Asymptotic Behavior of Best approximations of Classes of Periodic Analitic Functions Defined by Moduli of Continuity //  Bulgarian-Turkish-Ukrainian Scientific Conference $"$Mathematical Analysis, Differential Equations and their Applications\,$"$, Sunny Beach 15-20 September, 2010. -- Sofia: Academic Publishing House $"$Prof. Marin Drinov$"$, 2011. -- P.~173--182.

\item \label{Serd 2005} {\it Сердюк А.С.} Наближення класів аналітичних функцій сумами Фур'є в рівномірній метриці  // Укр. мат. журн. -- 2005. -- {\bf 57}, №8.~-- C.~1079--1096.

\item \label{Serduk 2012} {\it Сердюк А.С.} Наближення інтерполяційними тригонометричними поліномами на класах періодичних аналітичних функцій  // Укр. мат. журн. -- 2012.~-- {\bf 64}, №5.~-- C.~698--712.

\item \label{Gren} {\it Грандштейн И.С., Рыжик И.М.} Таблицы интегралов, сумм, рядов и  произведений. -- М.: Наука, 1971. --  1108 с.

\item \label{Serduk Sokolenko 2011}{\it  Сердюк А.С., Соколенко І.В.\/} Рівномірні наближення класів $(\psi,\overline{\beta})$-диференційовних функцій лінійними методами //  Теорія наближення функцій та суміжні питання: Зб. праць Ін-ту математики НАН України. -- 2011. -- 8, №1. -- С.~181--189.

\item \label{Rykasov 2003} {\it Рукасов В.И.} Приближение суммами Валле Пуссена классов аналитических функций //
Укр. мат. журн. -- 2003. -- {\bf 55}, №6.~--  \ C.~806--816.

\item \label{Korn} {\it Корнейчук Н.П.} Точные константы в теории приближения. -- М.: Наука, 1987.~-- 422 c.

\item \label{Teliacovs 1989}{\it  Теляковский С.А.\/} О приближении суммами Фурье функций высокой гладкости// Укр. мат. журн.-- 1989.-- {\bf 41}, №~4.-- C.~510--518.

\end{enumerate}

\end{document}